\documentclass [11pt] {amsart}
\usepackage{amsthm,amsmath,amssymb,amscd,graphics,enumerate}
\usepackage{latexsym}
\usepackage{verbatim}
\usepackage{enumerate}
\usepackage{amsthm}
\usepackage{amsmath}
\usepackage{amscd}

\newcounter{theorem}[section]
\numberwithin{equation}{section}

\newtheorem{cl}[theorem]{Claim}

\newtheorem{Thm}[theorem]{Theorem}
{\theoremstyle{remark}
 
\newtheorem{Rem}[theorem]{\text{\textbf{Remark}}} }
\newtheorem{Def}[theorem]{Definition}
\newtheorem{Lem}[theorem]{Lemma}
\newtheorem{Prop}[theorem]{Proposition}
\newtheorem{Cor}[theorem]{Corollary}
\newtheorem{situation}[theorem]{Situation}
\newtheorem{notation}[theorem]{Notation}
 
{\theoremstyle{definition}
 }

\newcommand{\RR}{{\mathbb{R}}}

\newcommand{\QQ}{{\mathbb{Q}}}
\newcommand{\NN}{{\mathbb{N}}}

\newcommand{\ZZ}{{\mathbb{Z}}}

\newcommand{\bK}{{\mathbf{K}}}

\newcommand{\m}{{\mathfrak{m}}}

\newcommand{\Spec}{\operatorname{Spec}}

\newcommand{\mExt}{\mathcal{E}xt}
\newcommand{\mHom}{\mathcal{H}om}

\newcommand{\mO}{\mathcal{O}}  
\newcommand{\mX}{\mathcal{X}}
\newcommand{\mY}{\mathcal{Y}}

\newcommand{\mS}{\mathcal{S}}

\newcommand{\mC}{\mathcal{C}}

\newcommand{\mK}{\mathcal{K}}

\newcommand{\mM}{\mathcal{M}}
\newcommand{\mR}{\mathcal{R}}

\newcommand{\mZ}{\mathcal{Z}}
\newcommand{\fbar}{\bar{f}}

\title{Cone Theorem via Deligne-Mumford stacks}
\author{Jiun-Cheng Chen}
\address{Department of Mathematics\\ Northwestern University\\ 2033 Sheridan Road\\ Evanston\\ IL 60208-2370\\ USA}
\email{jcchen@math.northwestern.edu}

\author{Hsian-Hua Tseng}
\address{Department of Mathematics\\ University of California\\ Berkeley\\ CA 94720\\ USA}
\email{hhtseng@math.berkeley.edu}

\date{\today}

\begin{document}
\begin{abstract}
We prove the cone theorem for varieties with LCIQ singularities using deformation theory of stable maps into Deligne-Mumford stacks. We also obtain a sharper bound on $-(K_X+D)$-degree of $(K_X+D)$-negative extremal rays for  projective $\QQ$-factorial log terminal threefold pair $(X,D)$.
\end{abstract}
\maketitle
\section{Introduction}
We work over an algebraically closed field $k$ of characteristic $ch(k) \geq 0$. 

The Cone Theorem is an important theorem in birational geometry. For smooth varieties, it can be stated as follows:
\begin{Thm}[Mori \cite{mo82}]
Let $X$ be a smooth projective variety. Then  
\[ \overline{NE}(X)= \overline{NE}(X)_{K_{X} \geq 0}+ \sum \RR _{\geq 0}[C_i]\] 
for a countable collection of  curves $C_i$ (extremal rays) with $K_X \cdot C_i <0$. The cone   $\overline{NE}(X)$ is locally finite  in the open half-space $N_1(X)_{K_X <0}$, and $$-(dim X +1)\leq C_{i} \cdot (K_X) < 0.$$
\end{Thm}

 Mori's proof for smooth varieties used a geometric argument, the famous ``bend and break'' method. The smoothness of $X$ is essential in his arguments: it is used to control the deformations of maps from curves to $X$.  When $X$ is log terminal (and $ch (k)=0$), a cohomological approach has been developed by Kawamata, Reid, Shokurov, and Koll\'ar (see \cite{ka84}, \cite{r83}, \cite{sh85} and \cite{ko84}). This method yields the bound $2\; dim\;X$ on $-K_X$-degree of extremal rays, instead of $dim \;X+1$ (see Keel \cite{ke99} for some results in $ch(k)=p$ and $dim X=3$). This class of singularities is large enough for the purpose of the minimal model program (MMP). These two approaches are very different in nature. 
Koll\'ar \cite{ko92} later developed a deformation theoretical argument for the case when $X$ is projective and LCIQ (quotients of LCI singularities). For $X$ projective and LCIQ, he constructed a LCI algebraic space $X^b$ which he called the bug-eyed cover. He proved the cone theorem by studying deformations of maps from curves to $X^b$. His method also yields the bound $2\; dim \;X$ on $-K_X$-degree. He also obtained a very general form (more general than the cohomological approach) of the  Cone theorem in dimension $3$, which  required strong  results from  the MMP in dimension $3$. 
 
In this paper, we develop a deformation theoretical method  using Deligne-Mumford stacks, instead of bug-eyed covers. For $X$ projective and LCIQ, we consider the associated LCI Deligne-Mumford stack $\mX$. We prove the cone theorem by studying deformations of stable maps to $\mX$. One advantage of using Deligne-Mumford stack is that the space of stable maps into a proper Deligne-Mumford stack is known to be proper \cite{av02}. Philosophically, we can then take the limit when we degenerate curves. When $X$ has tame LCIQ singularities, the bound on $-K_{X}$-degree obtained by our method is $dim\;X+1$, the same as that in smooth case.   

More precisely, we obtain:
\begin{Thm}[=Theorem~\ref{cone1-0}]\label{1cone1-0}
Let $X$ be a projective variety with only tame LCIQ singularities.  Then there is a countable collection  ${\bf \mR}$ of   $K_{X}$-negative extremal rays such that
\[ \overline{NE}(X)=  \overline{NE}(X)_{K_{X} \geq 0}+ Im[\overline{NE}(X_{sing}) \to \overline{NE}(X)] + \sum_{R_{i}  \in  \mR} R_i.\]
Furthermore, for each ray $R_i \in {\bf \mR}$, there is a rational curve $[C_i] \in R_i $  such that  $$-2\; dim X \leq C_i \cdot K_X <0.$$  

\end{Thm} 
\begin{Rem}
A stronger version of Theorem~\ref{1cone1-0} has been obtained in \cite{ko92}.
\end{Rem}

\begin{Thm}[=Theorem~\ref{cone1}]\label{1cone1}
 Let $X$ be a projective variety with only tame LCIQ singularities.  Then there is a countable collection  ${\bf \mR}$ of $K_{X}$-negative extremal rays such that
\[ \overline{NE}(X)=  \overline{NE}(X)_{K_{X} \geq 0}+ Im[\overline{NE}(X_{sing}) \to \overline{NE}(X)] + \sum_{R_{i}  \in  \mR} R_i.\]
Furthermore, for each ray $R_i \in {\bf \mR}$, there is a rational curve $[C_i] \in R_i $  such that  $$-(dim X+1) \leq C_i \cdot K_X <0.$$  

\end{Thm} 

\begin{Thm}[=Theorem~\ref{cone2}]\label{1cone2}
%%%change label here
Let $X$ be a projective variety with only isolated LCIQ singularities. Then the set ${\bf \mR}$ of all $K_{X}$-negative extremal rays of $\overline{NE}(X)$ is countable and 
\[ \overline{NE}(X)=  \overline{NE}(X)_{K_{X} \geq 0} + \sum_{R_{i} \in  \mR} R_{i}.\]
These rays are locally finite in the half-space $N_{1}(X)_{K_{X} <0}$, and for each $K_X$-negative extremal ray $R_i$ there is a rational curve $l_i$ generating $R_i$ such that $l_i \cdot (-K_{X}) \leq dimX+1$.
\end{Thm} 
 
We remark that in \cite{ma02}, Remark-Question 10-3-6, the bound $dimX+1$ is refered as ``conjecturally sharp'' for $\QQ$-factorial terminal varieties.

Using the MMP in dimension $3$, we are able to prove the following  theorem: 
\begin{Thm}[=Theorem~\ref{cone2-3}]\label{cone2-3-0}
Let $(X,D)$ be a projective threefold pair with $D$ a boundary divisor. Assume that $(X,D)$ has  only divisorial  log terminal  singularities.
 Then the set ${\bf \mR}$ of all $(K_{X}+D)$-negative extremal rays of  $\overline{NE}(X)$ is countable and 
\[ \overline{NE}(X)=  \overline{NE}(X)_{(K_{X}+D) \geq 0} + \sum_{R_{i} \in  \mR} R_{i}.\]
These rays are locally finite in the half-space $N_{1}(X)_{(K_{X}+D) <0}$, and for each ray $R_i$ there is a rational curve $l_i$ generating $R_i$ such that $$-4 \leq l_i \cdot (K_{X}+D) <0.$$
\end{Thm}
For a general $\QQ$-factorial projective threefold  $X$, we do not expect such a bound.

\subsection{Setting} 

\begin{Def}
\hfill
\begin{enumerate}
\item A scheme $Y$  is LCI if for every $y \in Y$, there is an \'etale neighborhood $U$ of $y$ such that $U$ is a complete intersection. This notion also makes sense when $Y$ is a Deligne-Mumford stack. 
\item We say that $Y$ has only LCIQ singularities if for every $y \in Y$, there is an \'etale neighborhood $U$ of $y$  such that $U \cong V/G_y$ where $V$ is LCI 
and the action of $G_y$ on $V$ is \'etale in codimension $1$.  
\item Notation as in (2). We say $Y$ has only tame LCIQ singularities if either $ch(k)=0$ or $ch(k)=p>0$ and the order $\mid G_y \mid$ of $G_y$ is not divisible by $p$.  
\end{enumerate}
\end{Def}
\begin{Def}[The Mori cone $\overline{NE}(Y)$]
Let $Y$ be a projective variety. Let $N_1(Y)$ denote the finite dimensional vector space 
$$ \{1-\text{cycles with real coefficients}\}/ \{\text{numerical equivalence}\}.$$ Let $NS(Y)$ be the N\'eron-Severi group of $Y$. The vector space $NS(Y) \otimes \RR$ is dual to $N_1(Y)$. Consider the cone $NE(Y) \subset N_1 (Y)$ generated by the numerical equivalence classes of effective $1$-cycles. Denote by $\overline{NE}(Y)$ the closure of $NE(Y)$ in $N_1 (Y)$.
\end{Def}
\begin{Def} 
Let $v \in \overline{NE}(Y) \subset N_1 (Y)$. The ray $R= \RR_{\geq 0} v \subset \overline{NE}(Y) $ is called extremal if 
\[ u_1, u_2 \in \overline{NE}(Y),  u_1+ u_2 \in R \Rightarrow u_1, u_2 \in R. \]\end{Def}

\begin{comment}
If an extremal ray $R = \RR _{\geq 0} [C]$ where $C$ is a curve on $Y$. 
We will choose a curve $C'$ such that $R=  \RR _{\geq 0} [C']$, and $C'$ has the minimal degree with respect to a given ample divisor on $Y$.
 
In this paper, $X$ is normal, projective and LCIQ. In some  cases, we also assume it has only isolated singularities. The  Deligne-Mumford stack $\mX$ is assumed to have a projective coarse moduli space $X$.   
% need to write LCI, LCIQ
\end{comment}

\section*{Acknowledgments}
We thank Dan Abramovich, Kai Behrend, Matthew Emerton, Tom Graber, J\'anos Koll\'ar, Jun Li, Ravi Vakil and Angelo Vistoli for helpful discussions. Part of this work was done when the authors visited National Center of Theoretical Sciences in Hsinchu, Taiwan. It is a pleasure to acknowledge their hospitality and support.

\section{Twisted stable maps}
Let $\mX$ be a proper tame  Deligne-Mumford stack with projective coarse moduli space $X$. The notion of twisted stable maps to $\mX$ has been studied by Abramovich and Vistoli \cite{av02}. In this section, we summarize the definitions and properties of twisted stable maps and their moduli  for the reader's convenience.
\begin{notation}
Unless otherwise mentioned, $\pi: \mX \to X$ denotes the map to the coarse moduli space.
\end{notation}

\subsection{Twisted curves}

An important ingredient of the paper \cite{av02} is what they called {\em twisted curves} (we sometimes call these orbicurves). Roughly speaking, these are nodal curves having certain stack structures \'etale locally near nodes (and, for pointed curves, marked points). We describe these stack structures in more details below. For the precise definition, see \cite{av02}, Definition 4.1.2. We restrict ourselves to {\em balanced} twisted curves throughout the paper. Let $\mC$ be a twisted curve and $C$ its coarse moduli space.

\subsubsection{Nodes}
For a positive integer $r$, let $\mu_r$ denote the cyclic group of $r$-th roots of unity. \'Etale locally near a node, a twisted curve $\mC$ is isomorphic to the stack quotient $[U/\mu_r]$ of the nodal curve $U=\{xy=f(t)\}$ by the following action of $\mu_r$: $$(x,y)\mapsto (\zeta_r x,\zeta_r^{-1}),$$ where $\zeta_r$ is a primitive $r$-th root of unity. \'Etale locally near this node, the coarse curve $C$ is isomorphic to the schematic quotient $U/\mu_r$.

\subsubsection{Markings}
\'Etale locally near a marked point, $\mC$ is isomorphic to the stack quotient $[U/\mu_r]$. Here $U$ is a smooth curve with local coordinate $z$ defining the marked point, and the $\mu_r$-action is defined by $$z\mapsto \zeta_r z.$$ Near this marked point the coarse curve is the schematic quotient $U/\mu_r$. 
\subsection{Twisted stable maps}
\begin{Def}
A twisted $n$-pointed stable map of genus $g$ and degree $d$ over a scheme $S$ consists of the following data (see \cite{av02}, Definition 4.3.1):
$$\begin{CD}
\mC @>f>>  \mX \\
@V{\pi_\mC}VV @V{\pi}VV \\
C @>\fbar>> X \\
@V{}VV \\
S.
\end{CD}$$
\begin{comment}
$$ \begin{array}{ccc}
 \mC             & \overset{f}{\to}&  \mX \\
\pi_\mC \downarrow     &    & \pi \downarrow \\
 C             & \overset{\fbar}{\to}& X \\
\downarrow     &    &\\
S \end{array}$$
\end{comment}
along with $n$ closed substacks $\Sigma_i\subset \mC$ such that
\begin{enumerate}
\item $\mC$ is a twisted nodal $n$-pointed curve over $S$ (see \cite{av02}, Definition 4.1.2),
\item $f:\mC\to \mX$ is representable,
\item $\Sigma_i$ is an \'etale gerbe over $S$, for $i=1,...,n$, and
\item the map $\fbar: (C,\{p_i\})\to X$ between coarse moduli spaces induced from $f$ is a stable $n$-pointed map of degree $d$ in the usual sense (see for instance \cite{FuP}).
\end{enumerate}
\end{Def}
A word on stability: $f:\mC\to \mX$ is stable if and only if for every irreducible component $\mC_i\subset \mC$, one of the following holds:
\begin{enumerate}
\item
$f|_{\mC_i}$ is nonconstant,
\item
$f|_{\mC_i}$ is constant, and $\mC_i$ is of genus at least $2$,
\item
$f|_{\mC_i}$ is constant, $\mC_i$ is of genus $1$, and there is at least one special points on $\mC_i$,
\item
$f|_{\mC_i}$ is constant, $\mC_i$ is of genus $0$, and there are at least three special points on $\mC_i$.
\end{enumerate}
In particular, a nonconstant representable morphism from a smooth twisted curve to $\mX$ is stable.

We say a twisted stable map $\mC \to \mX$ is rational if the coarse moduli space $C$ of $\mC$ is rational.

Let $\mK_{g,n}(\mX,d)$ denote the category of twisted $n$-pointed stable maps to $\mX$ of genus $g$ and degree $d$. The main result of \cite{av02} is that $\mK_{g,n}(\mX,d)$ is a proper Deligne-Mumford stack with projective coarse moduli space denoted by $\bK_{g,n}(\mX,d)$.

\begin{Rem}
Consider a line bundle $L$ on $X$. We have $\mC \cdot (f\circ \pi)^{*} L= \mC \cdot (\pi_{\mC} \circ \fbar) ^{*} L$, which is the same as $C \cdot  \fbar ^{*} L$ by the projection formula.
\end{Rem}

\subsection{Quotient singularities and Deligne-Mumford stacks}
\begin{Prop}\label{rightstack}
Let $X$ be a normal algebraic variety with LCIQ singularities. Then there is an LCI Deligne-Mumford stack $\mX$ such that $X$ and $\mX$ are isomorphic in codimension $1$ and $X$ is a coarse moduli space of $\mX$.
\end{Prop}
\begin{proof}
This is essentially \cite{vi89}, Proposition 2.8. We reproduce the arguments here. 

By assumption, there is a finite collection of varieties $V_i$ and morphisms $V_i\to X$ such that 
\begin{enumerate}
\item
$V_i$ is LCI,
\item
the morphisms $V_i\to X$ are \'etale in codimension $1$,
\item
for each $i$, there is a finite group $G_i$ acting on $V_i$ such that the morphism $V_i\to X$ factors as $V_i\to V_i/G_i\to X$ with $V_i/G_i\to X$ \'etale.
\item 
the images of $V_i\to X$ cover $X$.
\end{enumerate}
Let $V_{ij}$ be the normalization of $V_i\times_X V_j$. The data 
$$\coprod_{i,j} V_{ij}\rightrightarrows \coprod_i V_i$$
defines an \'etale algebraic groupoid. Let $\mX$ be its associated Deligne-Mumford stack. By construction, $X$ and $\mX$ are isomorphic in codimension $1$. By \cite{gi84}, Proposition 9.2, $X$ is a coarse moduli space of $\mX$. 
\end{proof}
\begin{Rem}
If $X$ has only tame LCIQ singularities, then the stack $\mX$ constructed in Proposition \ref{rightstack} is also tame.
\end{Rem}

\begin{Lem}\label{candiv}
In the situation of Proposition \ref{rightstack}, let $\pi:\mX\to X$ be the projection to the coarse moduli space. Then $\pi^{*} K_{X} \cong K_{\mX}$ as $\QQ$-Cartier divisors.
 \end{Lem} 
\begin{proof} 
Observe that $X$ and $\mX$ are isomorphic in codimension  $1$. Denote by $U \subset X$ the maximal open subset on which $X$ and $\mX$ are isomorphic. Since $K_{X}$ is $\QQ$-Cartier, take $m \in \mathbb{N}$ such that $mK_{X}$ is Cartier. Consider the pull back $\pi^*(mK_{X})$ to $\mX$. Since $X$ and $\mX$ are isomorphic over $U$,  $\pi^*(mK_{X})$ is isomorphic to $m K_{\mX}$ over $U$. It then follows that  $\pi^*(mK_{X}) \cong m K_{\mX} $ since $Codim(X\setminus U) \geq 2$. 
\end{proof}

\subsection{Lifting}
Let $X$ be a normal projective variety with LCIQ singularities and $\mX$ the stack as in Proposition \ref{rightstack}. Let $C$ be a smooth curve and $\fbar: C \rightarrow X$ a morphism whose image intersects the smooth locus of $X$. It follows that only finitely many points on $C$ are mapped to the singular locus of $X$. We want to ``lift''  the map $\fbar:C \to X$ to a map $C \to \mX$. It is not possible in general. However, we can always find a lifting after endowing a orbicurve structure on $C$. More precisely, let $\{p_i\}\subset C$ be the set of points which are mapped to the singular locus of $X$, and let $C_0=C\setminus \{p_i\}$. Since $X$ is isomorphic to $\mX$ outside the singular locus $X_{sing}$, the map $\fbar|_{C_0}:C_0\to X$ admits a lifting $C_0\to \mX$. By \cite{av02}, Lemma 7.2.5, there exists a twisted curve $\mC$ with coarse moduli space $C$, and a twisted stable map $f:\mC\to \mX$ extending $C_0\to \mX$. 

\section{Bend-and-Break}
We recall some results on Bend and Break for the reader's convenience. 
\begin{Prop}[see \cite{de} or \cite{km98}]\label{bandb1}
Let $X$ be a projective variety, $f: C \rightarrow X$ a smooth curve and $c$ a point on $C$. Assume $C$ is irrational.
If $dim_{[f]}Mor(C,X;f \mid_{c}) \geq 1$, then 
there is a curve $g: \tilde{C} \to X$ such that
\begin{enumerate}
\item the curve $\tilde{C}$ has  at least one rational component $\tilde{C}_1$ such that  $g \mid_{\tilde{C}_1}: \tilde{C}_1 \to X$ is non-constant and $f(c) \in  g_* \tilde{C}_1$, and 
\item $f_{*}C \sim g_* \tilde{C}$.
\end{enumerate}
\end{Prop}
\begin{Prop}[see \cite{de} or \cite{km98}]\label{bandb2}
Let $X$ be a projective variety, and $f: P^{1} \rightarrow X$ a rational curve. If $dim_{[f]}Mor(P^{1},X;f \mid_{\{0, \infty\}}) \geq 2$, then the $1$-cycle $f_{*}P^{1}$ is (algebraically) equivalent to a connected non-integral effective rational $1$-cycle passing through $0$ and $\infty$.
\end{Prop}

\begin{Prop}[see \cite{de}]\label{bandb3}
Let $X$ be a projective variety, $H$ a nef divisor, $f: C \rightarrow X$ a smooth curve, $B$ a finite subset  of $C$.  If $dim_{[f]}Mor(C, X,f \mid_{B} ) \geq 1$, then there is a rational curve $R$ which meets  $f(B)$, and  $H \cdot R \leq \frac{2H \cdot C}{Card(B)}$.   
\end{Prop}

\begin{Rem}
These results do not dependent on the smoothness of $X$. However, certain smoothness are needed to verify the condition that $C$ has enough deformations. 
\end{Rem}

We can extend these results to a proper Deligne-Mumford stack $\mX$ with projective coarse moduli space $X$. 
\begin{comment}
We assume that $\pi: \mX \to X$ is isomorphic in codimension $1$. When $ch(k)>0$, we assume that $\mX
$ is tame.
\end{comment}

\begin{Prop}\label{bandb2-0}
Let $\mX$ be a tame proper Deligne-Mumford stack with a projective coarse moduli space $X$. Let $C$ be a smooth irrational curve, $c $ a point on $C$ and $f: (C, c) \rightarrow \mX$ a twisted stable map. If $dim_{[f]}Mor(C,\mX;f \mid_{c}) \geq 1$, then there is a twisted stable map $f_1: \tilde{C} \to \mX$ such that 
\begin{enumerate}
\item 
the source curve $\tilde{C}$  has at least one rational component $\tilde{C}_1$ which has  at most two stacky points such that   $f_1 \mid_{\tilde{C}_1}: \tilde{C}_1 \to \mX$ is non-constant and $f(c) \in  f_{1*}\tilde{C}_1$, and
\item 
$f_{*}C \sim f_{1*}\tilde{C}$.
\end{enumerate}
\end{Prop}

\begin{Prop}\label{bandb2-1}
Let $\mX$ be a tame proper Deligne-Mumford stack with projective coarse moduli space $X$, $\mC$ an orbicurve with two marked points $0, \infty$ such that the coarse curve of $\mC$ is $P^1$ and $\mC$ has at most two stacky points. Let $f: \mC \rightarrow \mX$ be a $2$-pointed twisted stable map. If $dim_{[f]}Mor(\mC,\mX;f \mid_{\{0, \infty \}}) \geq 2$, then the $1$-cycle $f_{*}\mC$ is (algebraically) equivalent to a connected non-integral effective (twisted) rational $1$-cycle $\sum \mC_{i}$ passing through $0$ and $\infty$, 
 and one of its components has at most two stacky points.
\end{Prop}
\begin{comment}
\begin{Rem}
The conditions that $\pi: \mX \to X$ is isomorphic in codimension $1$ and $\mX$  is tame  are satisfied when 
 $X$ is covered by   \'etale neighborhoods $\{U_i\}$ such that $U_i \cong V_i/G_i$ where $V_i$
is  LCI, the action $G_i$ is \'etale in codimension $1$  and $p \nmid\;\mid G_i \mid$.
\end{Rem}
\end{comment}

 The proofs of Proposition~\ref{bandb2-0} and Proposition~\ref{bandb2-1} are very similar: The point is that we deform the map $f:\mC \to \mX$ with the source curve $\mC$ {\em fixed}. Any stable limit of such a (nontrivial) deformation is a map whose source curve is $\mC$ with a rational tree attached. We only prove Proposition~\ref{bandb2-1}.
\begin{proof}
Since $dim_{[f]}Mor(\mC,\mX;f \mid_{\{0, \infty\}}) \geq 2$, there is a nontrivial deformation   $F:\mC \times T \rightarrow \mX$  of $f$ over a smooth curve $T$. This induces a (non-constant) map $\phi: T\to \mK_{0,2}(\mX,d)$, where $d:=\pi_*f_*[\mC]$. Since $\mK_{0,2}(\mX,d)$ is proper, after a base change, there is a compactification $\bar{T}$ of $T$ and a map $\bar{\phi}: \bar{T}\to \mK_{0,2}(\mX,d)$ such that $\bar{\phi}|_T=\phi$. Let $\mS\to \mX$ be the $\bar{T}$-family of twisted stable maps corresponding to $\bar{\phi}$. Pick $t_0\in\bar{T}\setminus T$. Then the restriction to the fiber of $\mS\to \bar{T}$ over $t_0$ yields a $2$-pointed genus $0$ twisted stable map $f_{t_0}:\mC_{t_0}\to \mX$ (of nonzero degree). The domain curve $\mC_{t_0}$ is a tree of rational twisted curves. We claim that there is an irreducible component of $\mC_{t_0}$ such that 
\begin{enumerate}
\item \label{p1}
it has at most two special points,
\item \label{p2}
the restriction of $f$ to it is non-constant.
\end{enumerate} 
This can be seen as follows: Choose a maximal chain of rational twisted curves in $\mC$ and denote by $\mC', \mC''$ the two components at the ends of this chain. $\mC'\cup \mC''$ may contain both marked points, one of them, or none. If $\mC'\cup \mC''$ contains at most one marked point, then both $\mC'$ and $\mC''$ satisfy (\ref{p1}). If both $\mC'$ and $\mC''$ contain a marked point, then they both satisfy (\ref{p1}). If one of them, say $\mC'$, contains two marked points, then $\mC''$ contains no marked points and satisfy (\ref{p1}). In all cases, (\ref{p2}) follows from stability of $f$. 

The Proposition thus follows.
\end{proof}

\begin{Prop}\label{bandb4}
Let $\mC$ be a smooth orbicurve, $\Sigma\subset \mC$ the substack consisting of all stacky points, $\mX$ a tame proper Deligne-Mumford stack with projective coarse moduli space $X$, $\pi^* H$ a nef divisor on $\mX$ which is pulled back from a nef divisor $H$ on $X$, $f:\mC \rightarrow \mX$ a nonconstant proper representable morphism, and $B$ a finite subset of $C$ (may include some of the stacky points).  If $dim_{[f]}Mor(\mC, \mX, f|_B) \geq 1$, then $f_{*}\mC \sim \sum \mC_{i}$ where $\mC_{i}$'s are curves and at least one of them, say $\mC_{1}$, is a rational curve through a point on $f(B)$, and $\mC \cdot \pi^* H \leq \frac{2 C \cdot H}{Card(B)}$.     
\end{Prop}
\begin{proof}
Let $F:\mC \times T \rightarrow \mX$ be a nontrivial deformation of $f$ over a smooth curve $T$. This induces a non-constant map $\phi: T \to \mK_{g,b}(\mX,d)$ where $g$ is the genus of $\mC$, $b=Card (B\cup \Sigma)$ and $d:=\pi_*f_*[\mC]$. As in the proof of Proposition~\ref{bandb2-1}, after a base change we can compactify $T$ to a complete smooth curve $\bar{T}$ so that there is a morphism $\bar{\phi}:\bar{T}\to \mK_{g,b}(\mX,d)$ extending $\phi$. 
Denote the corresponding $\bar{T}$-family of twisted stable maps by $F: \mS \to \mX$, and the induced   $\bar{T}$-family of stable maps to $X$ by $\bar{F}: S \to X$. 
For any  $t_0\in \bar{T}\setminus T$,  let $\mS_{t_0}, S_{t_0}$ denote the fibers of $\mS\to \bar{T}, S\to \bar{T}$ over $t_0$ respectively. Clearly there is a one-to-one correspondence between the irreducible components of $\mS_{t_0}$ and those of $S_{t_0}$.
By \cite{ko96}  Theorem II-5.4, there is a rational component $C_1 \subset S_{t_0}$ such that $C_1 \cdot H \leq \frac{2 C \cdot H}{Card (B)}$ for at least one $t_0 \in  \bar{T} \setminus T$. Let $\mC_1$ be the corresponding component of $\mS_{t_0}$.
The proposition follows since  $C_1 \cdot H = \mC_1 \cdot \pi^* H$ and $C \cdot H= \mC \cdot \pi^*  H$.
\end{proof}

\begin{Rem}
The orbicurve $\mC_1$ may have more than two stacky points. 
\end{Rem}

\section{Deformations}
The following lemma is from \cite{ko92}, which is a generalization of the smooth case.
\begin{Lem}[\cite{ko92}]\label{def1}
Let $C$ be a proper connected algebraic curve without embedded points. Let $f: C \rightarrow Y$ be a morphism to an algebraic space $Y$ of pure dimension $n$. Assume that $Y$ is LCI, and the image of every component of 
$C$ intersects the smooth locus of $Y$. Then $dim_{[f]} Mor(C,Y) \geq -C \cdot K_{Y}+n \chi(\mO_{C})$.
\end{Lem}
\begin{Rem}
The condition that no component of $C$ lies completely in the singular locus is important. 
\end{Rem}

We extend this result to Deligne-Mumford stacks. 
\begin{situation}\label{curveset}
$\mC$ is a proper twisted curve, $\Sigma\subset \mC$ is the closed substack consisting of all smooth stacky points. $\mY$ is a proper Deligne-Mumford stack of pure dimension $n$ with projective coarse moduli space. $f:\mC\to \mY$ is a representable morphism such that the images of every component of $\mC$ under $f$ intersect the smooth locus of $\mY$. $B\subset \mC$ is a substack consisting of finitely many distinct smooth points (twisted or untwisted). $I_B$ is the ideal sheaf of $B$.
\end{situation}

By \cite{i71}, the space of first-order deformations of $[f]$ fixing $B$ is $Ext^0(f^*\Omega_\mY,I_{B})$, and the obstruction lies in $Ext^1(f^*\Omega_\mY,I_{B})$.

\begin{Lem}\label{def3}
In Situation \ref{curveset}, assume $\mY$ is smooth. Then 
$$dim_{[f]}Mor(\mC,\mY;f|_B)\geq n\chi(\mO_\mC)-\mC\cdot K_\mY-\sum_{x\in \Sigma\setminus B} age(f^*T\mY,x)-nCard(B).$$
\end{Lem}
\begin{proof}
By \cite{ka95}, the dimension $dim_{[f]}Mor(\mC,\mY;f|_B)$ is at least 
\begin{equation}\label{lowerbound}
dim Ext^0(f^*\Omega_\mY,I_{B})- dim Ext^1(f^*\Omega_\mY,I_{B}).
\end{equation}
For $\mY$ smooth, we have $Ext^i(f^*\Omega_\mY,I_{B})=H^i(\mC,f^*T\mY\otimes I_{B})$. Apply Lemma \ref{eulerchar} below to $V=f^*T\mY$ and observe that if $x\in \mC$ and $x\simeq B\mu_r$, then $age(V,x)+dim V|_x^{\mu_r}\leq dim V_x$.
It follows that (\ref{lowerbound}) is at least 
$$\chi(\mO_\mC)n-K_\mY\cdot \mC-\sum_{x\in \Sigma\setminus B} age(f^*T\mY,x)-nCard(B).$$
\end{proof}

\begin{Lem}\label{eulerchar}
Let $\mC$ and $B$ be as in Situation \ref{curveset}. Assume $B=\cup_i p_i$ with $p_i\simeq B\mu_{r_i}$. Let $V$ be a locally free sheaf on $\mC$. Then 
$$\chi (V\otimes I_{B})=\chi(\mO_\mC)\text{rank}V+\text{deg} V-\sum_{x\in \Sigma}\text{age}(V,x)-\sum_i \text{dim}V|_{p_i}^{\mu_{r_i}}.$$
\end{Lem}

\begin{proof}
Consider the exact sequence $0\to I_{B}\to \mO_\mC\to \mO_{B}\to 0$. Tensoring with $V$ and taking cohomology yield an exact sequence 
$$0\to H^0(\mC,V\otimes I_{B})\to H^0(\mC,V)\to H^0(\mC,V\otimes \mO_{B})$$
$$\to H^1(\mC,V\otimes I_{B})\to H^1(\mC,V)\to H^1(\mC,V\otimes \mO_{B}) \to 0.$$
Since $dim B=0$, we have $H^1(\mC,V\otimes \mO_{B})=0$. Taking alternating sum of dimensions, we obtain
$$dim H^0(\mC,V\otimes I_{B})-dim H^1(\mC,V\otimes I_{B})$$
$$= dim H^0(\mC,V)-dim H^1(\mC, V) -dim H^0(\mC,V\otimes \mO_{B}).$$
By Riemann-Roch for stacks (see \cite{Kw79}, \cite{t99}, and in particular \cite{ca04}, Lemma 2.4), we have
$$dim H^0(\mC,V)-dim H^1(\mC, V)=\chi(\mO_\mC)\text{rank}V+\text{deg} V-\sum_{x\in \Sigma}\text{age}(V,x).$$
Observe that $$H^0(\mC,V\otimes \mO_{B})=\oplus_i V|_{p_i}^{\mu_{r_i}}.$$ The Lemma is proved.
\end{proof}

\begin{Lem}\label{def2-1}
In Situation \ref{curveset}, assume $\mY$ is LCI. Then 
$$dim_{[f]}Mor(\mC,\mY;f|_B)\geq n\chi(\mO_\mC)-\mC\cdot K_\mY-nCard(\Sigma\cup B).$$
\end{Lem}

\begin{proof}
Applying the arguments of Theorem 2.10 of \cite{ko92}, we find that $f^*\Omega_\mY$ has projective dimension $1$. By the results of \cite{to02}, there is a global resolution of $f^*\Omega_\mY$ by locally free sheaves,
\begin{equation}\label{resolution}
0\to E\to F\to f^*\Omega_\mY\to 0.
\end{equation}
Applying $Hom_\mC (-,I_{B})$ and taking cohomology yield an exact sequence
$$0\to Hom(f^*\Omega_\mY,I_{B})\to Hom(F,I_{B})\to Hom(E,I_{B})$$
$$\to Ext^1(f^*\Omega_\mY,I_{B}) \to Ext^1(F,I_{B}) \to Ext^1(E,I_{B}) \to Ext^2(f^*\Omega_\mY,I_{B}).$$
We claim that $Ext^2(f^*\Omega_\mY,I_{B})=0$: Consider the local-to-global spectral sequence
$$E_2^{i,j}=H^i(\mC,\mExt^j(f^*\Omega_\mY,I_{B}))\Rightarrow Ext^{i+j}(f^*\Omega_\mY,I_{B}).$$
Note that $H^2(\mC, \mHom (f^{*} \Omega_\mY, I_{B}))=0$ since $dim\mC=1$ and $\mC$ is tame. $H^1(\mC, \mExt ^1( f^{*} \Omega_\mY, I_{B}))=0$ since $\mExt ^1 ( f^{*} \Omega_\mY, I_{B})$ is supported on isolated points (since the image of $f$ intersects the smooth locus of $\mY$). Also, since $f^{*} \Omega_\mY$ has projective dimension $1$, $\mExt ^2 ( f^{*} \Omega_\mY, I_{B})=0$. Hence $H^0(\mC, \mExt ^2 ( f^{*} \Omega_\mY, I_{B}))=0$. We conclude by the spectral sequence that $Ext^2(f^*\Omega_\mY,I_{B})=0$.

It follows that (\ref{lowerbound}) is equal to $\chi(F^\vee\otimes I_{B})-\chi(E^\vee\otimes I_{B})$. Applying Lemma \ref{eulerchar}, we find 
$$\chi(F^\vee\otimes I_{B})= \chi(\mO_\mC)rank F^\vee+deg F^\vee -\sum_{x\in\Sigma} age(F^\vee,x)-\sum_i dimF^\vee|_{p_i}^{\mu_{r_i}},$$
$$\chi(E^\vee\otimes I_{B'})=\chi(\mO_\mC)rank E^\vee+deg E^\vee -\sum_{x\in\Sigma} age(E^\vee,x)-\sum_i dimE^\vee|_{p_i}^{\mu_{r_i}}.$$
Note that $rankF^\vee-rank E^\vee =n$ and 
$$deg E^\vee-deg F^\vee=deg(detE^{-1}\otimes detF)=deg f^*\omega_\mY.$$
For $x\in\mC$ with $x\simeq B\mu_r$, taking stalks of the sequence (\ref{resolution}) yields a $\mu_r$-equivariant exact sequence,
$$0\to E_x\to F_x\to f^*\Omega_\mY|_x\to 0.$$ 
Write $E_x=\sum_{0\leq l<r} E_x^l$, $F_x=\sum_{0\leq l< r}F_x^l$ where $E_x^l$ and $F_x^l$ denote the $\mu_r$-eigen subspaces of $E_x$ and $F_x$ respectively with eigenvalue $\zeta_r^l$. Then we have $E_x^l\hookrightarrow F_x^l$. Note that the age terms are defined as $$age(E,x)=\sum_{0\leq l<r}\frac{l}{r}dim E_x^l,$$
$$age(F,x)=\sum_{0\leq l<r}\frac{l}{r}dim F_x^l.$$
It follows that
$$(age(F^\vee,x)+dim F^\vee|_x^{\mu_r})-(age(E^\vee,x)+dim E^\vee|_x^{\mu_r})$$
$$=\sum_{0\leq l<r}\frac{r-l}{r}dim F_x^l-\sum_{0\leq l<r}\frac{r-l}{r}dim E_x^l$$
$$\leq \sum_{0\leq l<r}(dim F_x^l-dim E_x^l)=rank F^\vee-rank E^\vee=n.$$
Therefore $$\chi(F^\vee\otimes I_{B})-\chi(E^\vee\otimes I_{B})\geq \chi(\mO_\mC)n-K_\mY\cdot \mC-nCard(\Sigma\cup B).$$
\end{proof}

Together with Bend and Break results, we have:

\begin{Lem}\label{bound}
Let $X$ be a normal projective variety with LCI singularities and $L$ an ample divisor on $X$. Let $f:C \rightarrow X$ be a smooth curve such that $f(C)$ meets the smooth locus of $X$ and $K_{X} \cdot C <0$. Given any point $x$ on $f(C)$, there exists a rational curve $R$ on $X$ through $x$ with  $L \cdot R \leq 2 dimX\frac{L \cdot C}{-K_{X} \cdot C}$.
\end{Lem}

%%more conditions when ch=p%%%
\begin{Lem}\label{bound1}
Let $X$ be a normal projective variety with (tame) LCIQ singularities and $L$ an ample divisor on $X$. Let $f: C \rightarrow X$ be a smooth curve such that $f(C)$ intersects the smooth locus of $X$ and $C \cdot K_{X} <0$. Given any point $x$ on $f(C)$, there exists a rational curve $R$ on $X$ through $x$ with  
$R \cdot  L  \leq 2 dimX\frac{C  \cdot  L}{C \cdot (-K_{X})}$.
\end{Lem}
\begin{proof} 
We only prove this lemma when $ch (k)=p>0$. Once we prove the result in $ch(k)=p$, the result in $ch (k)=0$ follows by a standard argument (see for example \cite{km98}, \cite{de}): Suppose that everything is defined over a finitely generated subring $\ZZ \subset R \subset k$. The only thing needed to be additionally careful with is to avoid those maximal ideals $\m$ that the characteristics of $R/\m$ divide the orders of stabilizer groups. These maximal ideals are contained in a finite union of subvarieties of $\Spec\;R$.

Let $\mX$ be the tame Deligne-Mumford stack as in Proposition \ref{rightstack}. Consider a lifting $\tilde{f}:\mC \to \mX $ of $f: C \to X$. Observe that we may assume $\mC$ is untwisted: Replace $\mC$ by a ramified cover $h:D \to \mC$ if necessary, and note that the bound $ \frac{D \cdot \pi^* L}{D \cdot (-K_{\mX})}$ equals to  $\frac{\mC  \cdot \pi^* L}{\mC \cdot (-K_{\mX})}$. Therefore we are reduced to the case where $C\to X$ factors through $\mX$. We abuse the notation and let $f$ denote the map $C \to \mX$. Consider the $m$-th Frobenius map $F_m: C \to C$. The degree (with respect to $\pi^* L$) of the composition $f_m:C \to C \to \mX$ is $p^m (\pi^* L \cdot_{f} C)$. Let $b_m=[\frac{- p^m (C \cdot_{f} K_X)-2}{n}]+1-g$ where $g$ is the genus of $C$ and $n=dim X$. By Lemma \ref{candiv} and \ref{def2-1}, the dimension at $f_m:C \to C \to \mX$ of the space of maps from $C$ to $\mX$ with $b_m$ points fixed is at least $2$. By Proposition~\ref{bandb4}, we can find a rational curve $\mC_m$ through a general point  $x$ on $f(C)$, and an effective cycle $\mZ_m$ such that $f_{m *}C \sim \mC_m+ \mZ_m$ and $\mC_m \cdot \pi^* L \leq   2\frac{p^m (C \cdot  L)}{b_m}$. 
Let $C_m:= \pi_* \mC_m$.
The limit of  $\frac{p^m (C \cdot L)}{b_m}$ is 
$n\frac{C  \cdot \pi^* L}{ C \cdot (-K_{X})}$. Also observe that $\mC_m \cdot \pi^* L= C_m \cdot L$. Note that $C_m \cdot  L$ is always an integer. Therefore, we can find some $m$ such that $C_m \cdot L \leq 2n\frac{L \cdot C}{-K_{X} \cdot C}$. 
\end{proof}

\section{Cone Theorem}
We start with the following classical result.
\begin{Thm}($ch (k) \geq 0$)\label{cone1-0}
Let $X$ be a projective variety with only tame LCIQ singularities. Then there is a countable collection ${\bf \mR}$ of $K_{X}$-negative extremal rays such that
\[ \overline{NE}(X)=  \overline{NE}(X)_{K_{X} \geq 0}+ Im[\overline{NE}(X_{sing}) \to \overline{NE}(X)] + \sum_{R_{i}  \in  \mR} R_i.\]
Furthermore, for each ray $R_i \in {\bf \mR}$, there is a rational curve $[C_i] \in R_i $  such that  $$-2\; dim X \leq C_i \cdot K_X <0.$$  

\end{Thm} 
\begin{proof}
 Let $\mX$ be the Deligne-Mumford stack as in Proposition \ref{rightstack}. Let $f: \mC \rightarrow \mX$ be a non-constant twisted stable map with $\mC$ smooth. Theorem~\ref{cone1-0} follows from Lemma~\ref{bound1} by an ingenious argument in \cite{ko92}, Theorem 3.3.
\end{proof}

\begin{Rem}\label{con}
%Theorem~\ref{cone1-0} follows from Lemma~\ref{bound1}. 
Koll\'ar \cite{ko96} proved a stronger version of Theorem~\ref{cone1-0}. In $ch (k)=p$, his proof does not need to assume $X$ has only  tame LCIQ  singularities.
\end{Rem}

We can improve the bound $2 dim X$:
\begin{Thm}\label{cone1}
Let $X$ be a projective variety with only tame LCIQ singularities. Then there is a countable collection ${\bf \mR}$ of $K_{X}$-negative extremal rays such that
\[ \overline{NE}(X)=  \overline{NE}(X)_{K_{X} \geq 0}+ Im[\overline{NE}(X_{sing}) \to \overline{NE}(X)] + \sum_{R_{i}  \in  \mR} R_i.\]
Furthermore, for each ray $R_i \in {\bf \mR}$, there is a rational curve $[C_i] \in R_i $  such that  $$-(dim X+1) \leq C_i \cdot K_X <0.$$  

\end{Thm} 

\begin{proof}
We only need to prove that the bound of $-K_X$-degree is $dimX+1$. 

First assume that $ch(k)=p>0$.
Let $R$ be a $K_X$-negative extremal ray which is not in $Im[\overline{NE}(X_{sing}) \to \overline{NE}(X)]$. Choose a rational curve $f: C \rightarrow X$ such that $R=\RR_{\geq0}[C]$. We may assume that for any rational curve $C'$ such that $R =\RR_{\geq0}[C']$, we have $-K_X \cdot C \leq -K_X \cdot C'$. Choose a (possibly ramified) cover $\rho : D  \rightarrow C$ so that the composition $D\to C\to X$ factors through $\mX$ and let $g$ denote the genus of $D$. Denote the map by $f_{D}:D \rightarrow \mX$. Note that $\rho_* D =d C$ for some positive integer $d$ and $D \cdot K_{\mX}=d C \cdot K_X$. It follows that $D \cdot K_{\mX} <0$ since $C \cdot K_X <0$. Note that $D$ may not be rational anymore. 

Consider the $m$-th Frobenius map $F_{m}:D \rightarrow D$. The degree of $f_{D} \circ F_{m}: D \rightarrow D \rightarrow \mX$ is $p^{m}(-K_{\mX} \cdot D)$. We may assume that $p^{m}(-K_{\mX} \cdot D)-g\; dimX >2$. Then by Proposition \ref{bandb2-0}  we can break $D$ as $D \sim  \sum \mC_{i}$ such that at least one of these curves, say $\mC_1$, is rational and has at most two stacky points. Denote the map by $f_1: \mC_1 \to \mX$. Since $[\pi \circ f_{D} \circ F_{m} (D)] \in R$ and $R$ is an extremal ray, it follows that the class $[\pi \circ f_{1}(\mC_1)]$ also generates $R$. The image $\pi \circ f_{1*} (\mC_1)$ does not lie in the singular locus $X_{sing}$ since $R$ is not in $Im[\overline{NE}(X_{sing}) \to \overline{NE}(X)]$. By Proposition ~\ref{bandb2-1}, as long as $\mC_1 \cdot K_{\mX}= \mC_1 \cdot \pi ^* K_X= C_1 \cdot K_X < -(1+dimX) $ (where $C_1$ is the coarse curve of $\mC_1$), we can further break $\mC_1$ into a union of rational curves one of whose irreducible components has at most two stacky points. Since $X$ is projective, this process will stop at some point. This concludes the proof when $ch (k)=p>0$.
%%%%%%%%%
%\end{proof}
%%%%
%
%\begin{proof}[Proof of the $ch(k)=0$ part of Theorem~\ref{cone1}]

Now assume that $ch(k)=0$. 
The proof of this case utilizes properties of extremal contractions. 

Let $C \subset X$ be a curve representing an extremal ray not contained in $Im[\overline{NE}(X_{sing}) \to \overline{NE}(X)]$. In particular, the class $[C]$ is not contained in 
$$Im[\overline{NE}(X_{- \infty}) \to \overline{NE}(X)],$$ where $X_{- \infty} \subset X_{sing}$ is the non log canonical locus. It is also clear that $C \nsubseteq X_{sing}$. By \cite{am03}, there is an extremal contraction $\phi: X \to Y$ which is projective and contracts only curves whose class is in  
the ray $\RR_{\geq 0}[C]$. 

The scheme $X_{sing}$ has only finitely many components. 
Let $X_{sing}^i$ be any component of $X_{sing}$, $Z^i_1 \to X^i_{sing}$ the normalization of $X_{sing}^i$, and $Z^i_2$ the normalization of $\phi (X^i_{sing})$. Set $Z_1:= \coprod_i Z_1^i$ and $Z_2:= \coprod_i Z_2^i$. We have the diagram:
$$\begin{CD}
Z_1 @>n>>  X_{sing} @>i>> X \\
@V{\phi \mid_{Z_1}}VV @V{\phi \mid_{X_{sing}}}VV  @V{\phi}VV \\
Z_2 @>n>> \phi(X_{sing}) @>i>> Y. 
\end{CD}$$
Note that $\phi \mid_{X_{sing}}: X_{sing} \to Y$ does not contract any curves since the only curve class contracted by $\phi: X \to Y$ is in the ray $\RR_{\geq 0}[C]$, which is not in $Im[\overline{NE}(X_{sing}) \to \overline{NE}(X)]$. Being projective and quasi-finite, we conclude that the morphism $\phi \mid_{Z_1}:Z_1 \to Z_2$ is finite. 

Let $k:=dim N_1(X)$.
Fix ample line bundles $H_i, \; i=1,2, \cdots ,k$ on $X$ which generate the vector space $N^1(X)$. Fix a positive integer $r$ such that $r K_X$ is a line bundle. Fix a positive number $\epsilon <<1$ such that $C \cdot (K_X+ \epsilon H_1) <0$.

Choose a finitely generated subring $\ZZ \subset R \subset k$ such that $X$, $f: P^1 \to C \subset X $, $X_{sing} \subset X$, $\phi(X_{sing}) \subset Y$, $Z_1 \to X$, $Z_2 \to Y$, the line bundles $H_i, H^{-1}_i \; i=1, 2, \cdots ,k$, and  the morphisms $\phi: X \to Y$ and  $\phi_{Z_1}: Z_1 \to Z_2$  are defined over $R$. Write $\phi_{R}: X_{R} \to Y_{R}$ and $\phi_{Z_1,\;R}:Z_{1, \; R} \to Z_{2, \; R}$ for the corresponding morphisms over $R$. Denote by $X_{\m}$, $Y_{\m}$ and $X_{sing , \; \m}$ the geometric fibers of the mod $\m$ reduction of $X_R$, $Y_R$ and $X_{sing ,\; R}$ respectively, and by $\phi_{\m}: X_{\m} \to Y_{\m}$  the mod $\m$ morphism of $\phi_R: X_R \to Y_R$. Denote the mod $\m$ reduction of $\phi_{Z_1, \; R}: Z_{1, \; R} \to Z_{2, \; R}$ by 
$\phi_{Z_1 , \;\m}: Z_{1, \; \m} \to Z_{2, \; \m}$.
We first show that for a dense subset of maximal ideals $\m$, the morphism 
$\phi_{Z_1 , \;\m}: Z_{1, \; \m} \to Z_{2, \; \m}$ is finite. 
In particular, this implies that the morphism $X_{sing, \m } \to \phi_{\m}(X_{sing,\; \m})$ is finite. 
%%%$[C_{\m}]$ is not contained in  
%%%$$Im[\overline{NE}(\overline{X}_{\m \;sing}) 
%%%\to \overline{NE}(\overline{X}_{\m})].$$ 
Since $X \setminus X_{sing}$ is smooth, $X_{\m} \setminus X_{sing , \; \m}$ is smooth for a dense subset of maximal ideals $\m$ (this implies that $X_{\m, sing}\subset X_{sing, \m}$). 
Since $\phi \mid_{Z_1}: Z_1 \to Z_2$ is a finite morphism, it follows that 
$\phi_{Z_1,\m}: Z_{1,\; \m} \to Z_{2, \; \m}$ is also a finite morphism for a dense subset of maximal ideals $\m$. 

%%%say more 

Let $\mM$ be the set of all maximal ideals of $R$ satisfying the following conditions:
\begin{enumerate}
\item the geometric fiber of X over $\m$ has only tame LCIQ singularities;
\item for $i=1, 2, ... k$, the line bundle $(H_i)_{\m}$ is ample (here $(H_i)_{\m}$ is the mod $\m$ reduction of $H_i$ on the geometric fiber $X_{\m}$ over $\m$);
\item $C_{\m} \cdot (H_i)_{\m} =C \cdot H_i$ for all $i$, where $C_\m$ is the mod $\m$ reduction of $C$;
\item $(-rK_X)_{\m} \equiv (-r K_{X_{\m}})$ where $(-rK_X)_{\m}$ is the mod $\m$ reduction of the line bundle $(-rK_X)$ on the geometric fiber $X_{\m}$;
\item the open variety $X_\m \setminus  X_{sing , \; \m}$ is smooth;
\item the morphism $\phi_{Z_1,\m}: Z_{1,\; \m} \to Z_{2, \; \m}$ is 
finite.
%%%$[C_{\m}]$ is not contained in 
%%%$Im[\overline{NE}(\overline{X}_{\m \;sing}) \to 
%%%\overline{NE}(\overline{X}_{\m})]$.
\end{enumerate}
Note that $\mM$ is dense in $\Spec\; R$.

Let $k (\m):=dim N_1(X_{\m})$.
Consider $[C_{\m}] \in \overline{NE}(X_{\m})$ and write it as a non-negative linear combination of at most $k(\m)$ extremal rays,
$$[C_{\m}]= b_{1, \m} w_{1, \m} + b_{2, \m} w_{2, \m}+ \cdots  +b_{k(\m), \m} w_{k(\m), \m},$$ see \cite{ko96} Lemma II-4.10.4. 
Since $\phi_{\m ,\;*}([C_{\m}])=0$, it follows that 
$\phi_{\m ,\;*}(w_{i, \m})=0$ when $b_{i, \m} \neq 0$. 
%%When $b_{i, \m}=0$, 
%%%we can take $ w_{i, \m}$ 
%%%to be any of the negative extremal ray represented by a curve. 
Therefore for those $i$ with $b_{i, \m}\neq 0$, the class $w_{i, \m}$ is not contained in 
$Im[\overline{NE}(X_{\m ,\;sing}) \to \overline{NE}(X_{\m})]$. This shows that the class $[C_\m]$ can be written as a linear combination (with positive coefficients) of at most $k(\m)$ extremal rays, each of which is not contained in $Im[\overline{NE}(X_{\m ,\;sing}) \to \overline{NE}(X_{\m})]$.
  
Since $C_{\m} \cdot (K_{X_{\m}}+ \epsilon (H_1)_{\m})<0$, at least one of the classes  
$\{w_{i, \m}: i=1, 2, \cdots , k(\m) \}$ lies in the half-space 
$\overline{NE}(X)_{(K_{X_{\m}}+ \epsilon (H_1)_{\m}) <0}$.
Combining all terms from the half space 
$\overline{NE}(X_\m)_{(K_{X_{\m}}+ \epsilon (H_1)_{\m}) \geq 0}$ into one term 
$u_\m$, we can write $$[C_m]= u_{ \m} + b_{2, \m} [C_{2, \m}]+ \cdots +b_{k(\m), \m} [C_{k(\m), \m}]$$ 
where $u_{\m} \in \overline{NE}(X_{\m})_{(K_{X_{\m}}+ \epsilon (H_1)_{\m}) \geq 0}$ and $[C_{i, \m}] \in \overline{NE}(X_{\m})_{ (K_{X_{\m}}+ \epsilon (H_1)_{\m}) <0}$. That is, we write $[C_{\m}]$ as a non-negative linear combination of $k(\m)$ vectors such that at most one of these vectors belongs to the (closed) half-space $\overline{NE}(X_{\m})_{(K_{X_{\m}}+ \epsilon (H_1)_{\m}) \geq 0}$ and the remaining $k(\m)-1$ vectors are extremal. (Note that some $b_{i, \m}$ may be $0$ since we may need fewer vectors in the linear combination.)
 
By the positive characteristic case of Theorem~\ref{cone1}, we may choose $C_{i, \m}$ so that $C_{i, \m} \cdot (-K_{X_{\m}}) \leq dimX +1,\; i=2, 3, \cdots ,k.$  
This implies immediately that $C_{i, \m} \cdot (H_1)_{\m}$ is bounded above by 
$\frac{dimX+1}{\epsilon}$ (note this number is independent of $\m$).
 
We identify $\overline{NE}(X)$ with its image under the embedding $\psi: \overline{NE}(X)\to \RR^k$ defined by $$[D]\mapsto (D\cdot H_1,...,D\cdot H_k).$$
For any $\m \in \mM$, consider the (geometric) fiber $X_m$ over $\m$  and the linear map 
$$\psi_{\m}: \overline{NE}(X_m) \to \RR^k$$
given by
$ \psi_{\m}([D])=(D \cdot (H_1)_{\m}, D \cdot (H_2)_{\m}, \cdots , D \cdot (H_k)_{\m})$. The image of $\psi_\m$ is a closed subcone in $\RR^k$ containing no lines (by Kleiman's ampleness criterion). 
This cone may not lie in $\overline{NE}(X) \subset \RR^k$.
Set $v_{\m}= \psi_{\m}([C_\m])$ and $v_{i,\m}=  \psi_{\m}([C_{i, \m}])$. Note that $v_\m$ does not depend on $\m$ by the choice of $\mM$. We denote this common vector by $v$. 
Consider the vector $$v= \psi_{\m}(u_{ \m})+b_{2, \m} v_{2, \m} + \cdots +b_{k(\m), \m} v_{k(\m), \m}$$ in $\RR^k$. We can express $v$ as a non-negative combination of at most $k$ vectors from $\{ \psi_{\m}(u_{\m}),v_{2, \m}, v_{3, \m}, \cdots , v_{k(\m), \m} \}$. We can lift this expression of $v$ to an expression of $[C_\m]$: Renaming $[C_{j,\m}], \; j= 2,3,  \cdots, k(\m)$ if necessary, we can write 
$$[C_{\m}] = \tilde{u}_{\m}+  a_{2, \m} [C_{2, \m}]+a_{3, \m} [C_{3, \m}] + \cdots +a_{k, \m} [C_{k, \m}] +w_{\m}$$ where $\tilde{u}_{\m}$ is either $u_{\m}$ or one of the vectors from 
$\{ [C_{2, \m}], [C_{3, \m}], \cdots, [C_{k(\m), \m }]\}$ (therefore 
$\tilde{u}_{\m} \in \overline{NE}(X_{\m})$), $a_{i,\m}\geq 0$, and $w_{\m}$ is a (not necessarily effective) linear combination of 
$$\{ u_{\m}, [C_{2, \m}], [C_{3, \m}], \cdots , [C_{k(\m), \m }] \}$$ such that $\psi_{\m}(w_{\m})=0$. Also, $$v=\psi_{\m}([C_{\m}])= \psi_{\m}(\tilde{u}_{\m}) + a_{2, \m} v_{2, m} + \cdots +a_{k, \m} v_{k, \m}.$$  

Consider the morphism scheme 
$$\pi_R: Mor_{1 \leq d \leq \frac{dim X+1}{\epsilon}}(P^1_R, X_R) \to \Spec \; R$$ where $d$ denotes the $H_1$-degree. The $R$-scheme $Mor_{1 \leq d \leq \frac{dim X+1}{\epsilon}}(P^1_R, X_R)$ is quasi-projective 
and has only finitely many components. Consider the decomposition
$$Mor_{1 \leq d \leq \frac{dim X+1}{\epsilon}}(P_R^1, X_R)= \coprod Mor^{(n_1,n_2, \cdots ,n_k)} $$ according to the intersection numbers $(n_1,n_2, \cdots, n_k)$ with $H_1, H_2, \cdots , H_k$. 
%%That is, if $[f] \in Mor^{(n_1,n_2, \cdots ,n_k)}$, then 
%%$$P^1 \cdot f^* H_1=n_1, P^1 \cdot f^* H_2=n_2, 
%%\cdots P^1 \cdot, f^* H_k=n_k.$$
We may view each $(n_1,n_2, \cdots, n_k)$ as a vector in $\RR^k$. Note that there are only finitely many such vectors. Let $S$ be the finite set of all these vectors, i.e. 
$$S:=\{(n_1,...,n_k)|\exists [f]\in Mor_{1 \leq d \leq \frac{dim X+1}{\epsilon}}(P^1_R, X_R) \text{ such that } P^1 \cdot f^* H_i=n_i \text{ for } i=1,...,k\}.$$
Put $S^{k-1}:=\{(v_1,v_2, \cdots ,v_{k-1}): v_i \in S\}$. Note that $v_{i,\m}=\psi_{\m}([C_{i,\m}])\in S$ by the choice of $C_{i, \m}$. 

For $t=(w_2, w_3, \cdots, w_k) \in S^{k-1}$, let $\mM_t$ be the set of all maximal ideals in $\mM$ so that for each $i=2, 3, ...,k$ we have $\psi_{\m}(C_{i, \m})= w_i$. Since $S$ is a finite set, for any infinite sequence of maximal ideals $\{ \m_{j} \}$, there is an infinite subsequence $ \{ \m_{j_l} \}$ such that for each $i=2, 3, \cdots ,k$, the sequence of vectors $ \{ v_{i, \m_{j_l}} \}\subset S$ is a constant sequence $\{ v_i, v_i, v_i,  \cdots \}$. This shows that $\coprod_{ t \in S^{k-1}} \mM_t$ is dense in $\Spec \;R$. Since $S^{k-1}$ is finite, $\mM_{t_0}$ is dense for at least one $t_0 \in S^{k-1}$. Write $t_0= (v_2, v_3, \cdots , v_k)$. Therefore, for any $\m \in \mM_{t_0}$ and $i=2, 3, \cdots , k$, we have $v_{i, \m}=\psi_{\m} ([C_{i, \m}])= v_i$. 

Fix an integer $i$ with $2\leq i \leq k$. Note that for $\m\in \mM_{t_0}$, the curve $C_{i,\m}$ gives a point in the fiber over $\m$ of the restriction of $\pi_R$ to $Mor^{v_i}$ (recall that $v_i\in S$). Since $\mM_{t_0}$ is dense, the fiber over the generic point is also non-empty; hence it has a geometric point which corresponds to a rational curve $f_i:P^1 \to C_i \subset X$. 

To conclude the proof, we need the next claim:
\begin{cl}\label{infseq}
There is an infinite sequence of maximal ideals  $\{ \m_k\}\subset \mM_{t_0}$ such that any infinite subsequence of $\{ \m_k \}$ is not contained in any finite union of subschemes of $\Spec \;R$ unless one of these subschemes is $\Spec\; R$. 
\end{cl}

Assuming Claim~\ref{infseq}, we continue the proof as follows:
Take a sequence $\{ \m_j: \m_j \in \mM_{t_0}\}$ as in Claim \ref{infseq}. Taking a subsequence if necessary, we may assume that for each $i=2, 3, \cdots k$ the sequence of vectors $\{ v_{i, \m_{j}}: j \in \NN \}$ is a constant sequence $\{v_i\}$. Note that $H_i$ is ample on $X_{\m_j}$. Therefore the  coordinates of the vectors $\psi_{\m}(\tilde{u}_{\m_j}),a_{2, \m_j} v_{2}, a_{3, \m_j}v_3, \cdots, a_{k, \m_j} v_{k}$ are all non-negative, and are at most the coordinates of $v$ since $v=\psi_{\m_j}(\tilde{u}_{\m_j})+a_{2,\m_j}v_2+ \cdots + a_{k, \m_j} v_k$. So the sequences $\{\psi_{\m}(\tilde{u}_{\m_j})\}$, $\{a_{2, \m_j} v_{2}\}$,..., $\{a_{k, \m_j} v_{k}\}$ are all bounded. Since any bounded infinite sequence in $\RR^k$ has a infinite subsequence which converges, there is a subsequence $\{ \m_{j_l} \}$ of $\{ \m_j \}$ such that 
$$\lim_{l \to \infty} a_{i, \m_{j_l}}= a_i\in \RR_{\geq 0}, \;\; i=2,3, \cdots ,k$$ and 
$$\lim_{l \to \infty} \psi_{\m_{j_l}}(\tilde{u}_{\m_{j_l}})=u\in \RR^k.$$
It follows that $v= u + a_2 v_2+ a_3 v_3 + \cdots +a_k v_k$. 

Recall that we identify $\overline{NE}(X)$ with a subcone of $\RR^k$ via the embedding $\psi$. If $u \in \overline{NE}(X)$, then since $v$ is extremal, it follows that 
$u,\; v_2, \; v_3, \; \cdots v_k \in \RR_{\geq 0}\; v$. Note that $C_{i,\m}\cdot (H_j)_\m=C_i\cdot H_j$, so each $v_i$ is represented by a rational curve $[C_i]$ with $-K_X$-degree at most $dim X+1$. (We have no control on the $-K_X$-degree of $u$.) This concludes the proof.  

Suppose that $u $ is not in the cone $\overline{NE}(X)$ (therefore it is non-zero!). By Kleiman's criterion for ampleness we can find an ample line bundle $H= c_1 H_1+ c_2 H_2 + \cdots +c_k H_k$ where $c_1,...,c_k \in \ZZ$ (possibly negative) such that $u \cdot H <0$ (view $H$ as a  linear function on $\RR^k$). This line bundle is defined over $R$. The locus of all maximal ideals such that $H$ is not ample on $X_{\m}$ is contained in a finite union of closed subschemes of $\Spec \;R$. Therefore, by Claim \ref{infseq} there are still infinitely many maximal ideals in the subsequence $\{\m_{j_l}\}$ such that $H$ is ample on $X_{\m_{j_l}}$. Since $\psi_{\m_{j_i}}(\tilde{u}_{\m_{j_l}}) \to u$ and 
$u \cdot H <0$, it follows that $\psi_{\m_{j_l}}(\tilde{u}_{\m_{j_l}}) \cdot H <0$ for $l$ large enough. This is impossible since $\tilde{u}_{\m_{j_l}} \in \overline{NE}(X_ {\m_{j_l}})$, and $H$ is ample on $X_{\m_{j_l}}$. Hence $u$ has to be in the cone $\overline{NE}(X)$.
\end{proof}
\begin{proof}[Proof of Claim \ref{infseq}]
Since $R$ is finitely generated over $\ZZ$. The ring $R$ is countable as a set. Since all ideals are finitely generated, the set of all ideals is also countable. List the set of corresponding subschemes as $\{V_1= \Spec \; R/I_1, V_2= \Spec \; R/I_2, V_3= \Spec \; R/I_3,  \cdots \}$.   
We construct a sequence as follows:
For each $i \in \NN $, pick a maximal ideal $\m_i \in \mM_{t_0}$ such that the corresponding point $\Spec \; R/\m_i$ is not in the union of first $i$ subschemes  $ V_1, V_2, V_3, \cdots ,V_i$. Such a $\m_i$ exists for every $i$ since $\mM_{t_0}$ is dense. It is not hard to verify this sequence has the declared property. 
\end{proof}

The following results are immediate.
\begin{Thm}\label{cone2}
Let $X$ be a projective variety with only isolated LCIQ singularities. Then the set ${\bf \mR}$ of all $K_{X}$-negative extremal rays of  $\overline{NE}(X)$ is countable and 
\[ \overline{NE}(X)=  \overline{NE}(X)_{K_{X} \geq 0} + \sum_{R_{i} \in  \mR} R_{i}.\]
These rays are locally finite in the half-space $N_{1}(X)_{K_{X} <0}$, and for each $K_X$-negative extremal ray $R_i$ there is a rational curve $l_i$ generating $R_i$ such that $l_i \cdot (-K_{X}) \leq dimX+1$.
\end{Thm} 

\begin{Thm}\label{cone2'}
Let $X$ be a projective variety with only quotient singularities. Then the set ${\bf \mR}$ of all $K_{X}$-negative extremal rays of $\overline{NE}(X)$ is countable and 
\[ \overline{NE}(X)=  \overline{NE}(X)_{K_{X} \geq 0} + \sum_{R_{i} \in  \mR} R_{i}.\]
These rays are locally finite in the half-space $N_{1}(X)_{K_{X} <0}$, and for each $K_X$-negative extremal ray $R_i$ there is a rational curve $l_i$ generating $R_i$ such that $l_i \cdot (-K_{X}) \leq dimX+1$.
\end{Thm}

\subsection{Threefolds}
The goal of this subsection is to prove:
\begin{Thm}\label{cone2-3}
Let $(X,D)$ be a projective threefold pair with $D$ a boundary divisor. Assume that $(X,D)$ has only divisorial log terminal singularities. Then the set ${\bf \mR}$ of all $(K_{X}+D)$-negative extremal rays of $\overline{NE}(X)$ is countable and 
\[ \overline{NE}(X)=  \overline{NE}(X)_{(K_{X}+D) \geq 0} + \sum_{R_{i} \in  \mR} R_{i}.\]
These rays are locally finite in the half-space $N_{1}(X)_{(K_{X}+D) <0}$, and for each ray $R_i$ there is a rational curve $l_i$ generating $R_i$ such that $$-4 \leq l_i \cdot (K_{X}+D) <0.$$
\end{Thm}
The part concerning the structure of $\overline{NE}(X)$ is already known, see for instance \cite{km98}. We will only prove the statement on the bound of $-(K_X+D)$-degree. 
The proof is divided into several steps:

{\bf Step 1:} The variety $X$ is terminal, $D$ is empty.
\begin{Cor}\label{cone2-1-0}
Let $X$ be a  projective threefold  with only terminal singularities. 
Then the set ${\bf \mR}$ of all $K_{X}$-negative extremal rays of  $\overline{NE}(X)$ is countable and 
\[ \overline{NE}(X)=  \overline{NE}(X)_{K_{X} \geq 0} + \sum_{R_{i} \in  \mR} R_i.\]
The rays $\{R_i\}$  are locally finite in the half-space $N_{1}(X)_{K_{X} <0}$,  and for each $K_X$-negative extremal ray $R_i$ there is a rational curve $l_i$ generating $R_i$ such that $l_i \cdot (-K_{X}) \leq 4$.
\end{Cor} 
\begin{proof}
%We only prove the part on the bound of $-K_X$-degree.
Since $X$ is terminal, it has only isolated hyper-quotient singularities. The result now follows immediately from Theorem~\ref{cone2}. 
\end{proof}

{\bf Step 2:} The variety $X$ is canonical, $D$ is empty.
\begin{Cor}\label{cone2-1}
Let $X$ be a  projective threefold  with only canonical singularities. Then the set ${\bf \mR}$ of all $K_{X}$-negative extremal rays of  $\overline{NE}(X)$ is countable and 
\[ \overline{NE}(X)=  \overline{NE}(X)_{K_{X} \geq 0} + \sum_{R_{i} \in  \mR} R_i.\]
The rays $\{R_i\}$  are locally finite in the half-space $N_{1}(X)_{K_{X} <0}$,  and for each $K_X$-negative extremal ray $R_i$ there is a rational curve $l_i$ generating $R_i$ such that $l_i \cdot (-K_{X}) \leq 4$.
\end{Cor} 
\begin{proof} 
Again, we only prove the part on the bound of $K_X$-degree. 
By the MMP, there exist a terminal threefold $Y$ and a birational morphism $g:Y \to X$ such that $K_{Y}$ is $g$-trivial (i.e. we can extract all exceptional divisors with discrepancy $0$. This is first proved by M. Reid, see \cite{km98} Theorem 6.23). %%add reference
 Let $R$ be a $K_X$-negative extremal ray. Choose a rational curve $C$ such that $R= \RR_{\geq 0}[C]$.

Consider $Z \subset Y$ such that $g_{*}Z=d C$ for some positive integer $d$ (when $C \nsubseteqq g(Exc(g))$, we can take $Z$ to be the proper transform of $C$ in $Y$ and  $d=1$).  
If $C \cdot K_{X} <-4$, then $Z \cdot K_{Y}=Z \cdot g^{*}K_{X}=d C \cdot K_{X} <-4d$. 

 By  Corollary~\ref{cone2-1-0} we can write $[Z]=u+ \sum a_i[Z_i]$ where $u \in \overline {NE}(Y)_{K_{Y} \geq 0}$ (it may be $0$), $Z_i$ is a rational curve with $-4 \leq Z_i \cdot K_Y <0$ and $a_i \in \RR_{\geq 0}$. Now $$g_* u +\sum  a_i g_{*}[Z_{i}]=  [g_{*}Z]= d [C].$$
Note that at least one  $Z_i$, say $Z_1$, is not contracted by $g$. 
Write $g_* Z_1= e_1 C_1$ where $C_1$ is the image of $Z_1$ in $X$ with the reduced structure and $e_1 \geq 1$.
%%% since $Z \cdots K_Y <0$.
Since 
$R= \RR_{\geq 0} [C]$ is an extremal ray, the curve $C_1$ also generates $R$, and    
\[e_1 C_1 \cdot K_X=g_{*}Z_{1} \cdot K_{X}=Z_{1} \cdot g^{*}K_{X}= Z_{1} \cdot K_{Y} \geq -4.\] 

It follows immediately that $ -4 \leq C_1 \cdot K_X$.
\end{proof}

{\bf Step 3:} The pair $(X,D)$ is canonical.
The goal of this step is to prove the following result. 
\begin{Prop}\label{cone2-1-1}
Let $X$ be a $\QQ$-factorial projective threefold, and $D=\sum a_{i}D_{i}$ a boundary divisor (i.e. $0 \leq a_i \leq 1$). Assume that $(X,D)$ has only canonical singularities. Then the set ${\bf \mR}$ of all $(K_{X}+D)$-negative extremal rays of $\overline{NE}(X)$ is countable and 
\[ \overline{NE}(X)=  \overline{NE}(X)_{K_{X}+D \geq 0} + \sum_{R_{i} \in  \mR}  R_i.\]
These rays are locally finite in the half-space $N_{1}(X)_{(K_{X}+D) <0}$, and for each ray $R_i$ there is a rational curve $l_i$ generating $R_i$ such that $$l_{i} \cdot (-(K_{X}+D)) \leq 3+1.$$
\end{Prop} 

We need several simple lemmas.
\begin{Lem}\label{lc-1}
Let $(S, D= \sum d_{i} D_{i})$ be a $\QQ$-factorial log canonical surface where $D$ is a boundary divisor, i.e. $0 \leq d_i \leq 1$. Let $R$ be any $(K_S+D)$-negative extremal ray and $[C] \in R$ an irreducible curve with minimal $-(K_S+D)$-degree, then $-3 \leq C \cdot (K_{S}+D)$. 
\end{Lem}

\begin{proof} 
First consider the case when $S$ is smooth. 
Suppose that 
$C \cdot (K_S+D) <-3$. If $C \cdot D_i \geq 0$ for all $i$, then 
$C \cdot K_S \leq  C \cdot (K_S+D) < -3$. 
By bend and break and the fact that $[C]$ is extremal, 
 it follows that  there is a rational curve $C'$ such that $[C'] \in R$ and $-3 \leq C' \cdot K_S < 0$. 
Since $[C']$ and $[C]$ are both in $R$, we have $[C']=a[C]$ for some positive rational number $a$ and $C' \cdot D_i = a C \cdot D_i \geq 0$. Thus  $C' \cdot (K_S+D) \geq -3$. So 
$$C'\cdot (-(K_S+D))\leq 3 < C\cdot (-(K_S+D)).$$
This contradicts to the choice of $C$.

If $C \cdot D_i <0$ for some $i$, then $C=D_i$. 
We have
$$D_i \cdot (K_S+D)= D_i \cdot (K_S+d_iD_i+ \sum_{j \neq i}d_j D_j) \geq 
D_i \cdot (K_S+d_i D_i) \geq D_i \cdot (K_S+D_i).$$
By adjunction formula, $D_i \cdot (K_S+D_i)=2 p_a(D_i)-2 \geq -2$. This gives a contradiction.
This proves the case when $S$ is smooth.

When $S$ is not smooth, consider a log resolution $(Y, D')$ of $(S,D)$ such that $D'$ is a snc divisor. Let $0 < \epsilon <<1$. 
Note that $(Y, (1 - \epsilon)D')$ is terminal for $0 < \epsilon $. 
Run the $(K_Y+ (1- \epsilon )D')$-MMP over $(S,(1- \epsilon)D)$. Let $(Y_1, (1-\epsilon)D'')$ be the final outcome. Denote the morphism by  $h: (Y_1, (1-\epsilon)D'') \to (S, (1- \epsilon)D)$.  
The pair $(Y_1, (1-\epsilon)D'')$
is still terminal (hence $Y_1$ is smooth) and $K_{Y_1}+(1- \epsilon)D''$ is $h$-nef. Write 
$$K_{Y_1}+(1- \epsilon)D'' = h^* (K_{S}+(1- \epsilon)D)+ \sum a_i E_i,$$ where $\{E_i\}$ is the collection of all irreducible curves contracted by $h$. 
Since $K_{Y_1}+(1- \epsilon)D''$ is $h$-nef, it follows that $a_i \leq 0$.

Let $C \subset S$ be an irreducible curve representing a $(K_S+D)$-negative 
extremal ray $R$. 
We take $C$ to have the lowest $-(K_S+D)$-degree.
Denote by $C'$ the proper transform of $C$ in $Y_1$. 
Note that $C'$ is not $h$-exceptional,  $(Y_1, (1- \epsilon)D''- \sum a_iE_i)$ is log canonical and  $(1- \epsilon)D''- \sum a_iE_i$ is effective.
We have 
$C' \cdot [ K_{Y_1}+(1- \epsilon)D''-\sum a_iE_i]= C' \cdot h^*(K_{S}+(1- \epsilon)D)$. 
We write $[C']$ as a sum of $(K_{Y_1}+(1- \epsilon)D''- \sum a_iE_i)$-negative extremal rays (which are represented by curves) and an element in  $\overline{NE}(Y_1)_{K_{Y_1}+(1- \epsilon)D'' \geq 0}$.
Write $[C']= \sum b_i[C_i]+ u^+$ in $\overline{NE}(Y_1)$ as in the proof of Step 2. 
%We may assume that 
%$C_i$ is not contracted by $h$.
Since $h_*([C'])=[C]$ and $[C]$ is extremal, we can write  $[h_*(C_i)]=e_i [C]$ and $h_* (u^+)=d [C]$. Again, at least one of $C_i$, say $C_1$, is not contracted by $h$. Note that $e_1\geq 1$ (otherwise $h_*(C_1)$ has smaller $-(K_S+D)$ degree). By the smooth case, we may choose $C_1$ so that 
$C_1 \cdot  [K_{Y_1}+(1- \epsilon)D'' - \sum a_iE_i] \geq -3$.  
We have \[e_1 C \cdot (K_S+(1- \epsilon)D)=h_*(C_1) \cdot (K_S+(1- \epsilon)D)\]\[=  C_1 \cdot  [K_{Y_1}+(1- \epsilon)D'' - \sum a_iE_i] \geq -3.\] 
It follows that $-3 \leq C_1 \cdot (K_S+(1-\epsilon)D)$.
\end{proof}

\begin{Lem}\label{sur1-1}
Let $(X,D=\sum d_i D_i)$ be a $\QQ$-factorial log canonical surface. If $C \subset X$ is a irreducible curve and $C \cdot (K_X+D)<-3$, then $C \cdot K_X < -3$.
%%change it to -3
\end{Lem}
\begin{proof}
It suffices
 to exclude the possibility that $C \cdot D <0$.
For simplicity we only consider the case when $X$ is smooth; the general case follows easily  from applying the smooth case to the minimal 
resolution of $(X,D)$ (more precisely, $(X, (1- \epsilon)D)$ as in Lemma~\ref{lc-1}).
If $C \cdot D <0$, then $C=D_i$ for some  $i$ and $$C\cdot (K_X+D)\geq D_i\cdot (K_X+D_i)=2p_a(D_i)-2\geq -2.$$ This  contradicts to the condition  
$C \cdot (K_X+D)<-3$.
\end{proof}

\begin{Lem}\label{sur1}
Let $(X,D= \sum d_i D_i)$ be a $\QQ$-factorial log canonical surface. If $C \subset X$ is a rational curve and $C \cdot (K_X+D) <-3$, then $C \sim \sum  C_{i}$ where $C_i$ is a rational curve, and $C_{i} \cdot (K_X+D) \geq -3$ ($C_i$ may equal to $C_j$).
\end{Lem}
\begin{proof}
Assume that $X$ is smooth. By Lemma \ref{sur1-1}, 
we have $C \cdot D \geq 0$ and $C\cdot K_X<-3$. By bend and break, we have $C \sim \sum  C_{i}$ with  $C_i$ a rational curve, and $C_{i} \cdot K_X \geq -3$. 
If $C_i \cdot D <0$, then 
 $C_i=D_j$ for some $j$ and $C_i\cdot (K_X+D)\geq C_i\cdot (K_X+D_j)\geq -2$ by adjunction. If $C \cdot D \geq 0$,  then $C_i\cdot (K_X+D)\geq C_i\cdot K_X\geq -3$.

If $X$ is not smooth, take a log resolution $(Y, D')$. We conclude by applying bend and break on $Y$ and using an argument similar to that in the proof of Lemma \ref{lc-1}.
\end{proof}

\begin{Lem}\label{sur2}
Let $X$ be a projective normal surface such that $K_X$ is $\QQ$-Cartier, and $Z$ an irreducible curve on $X$. Let $g:Y \to X$ be the minimal resolution of $X$ and $Z'$ the proper transform of $Z$. If $Z \cdot K_X <0$, then $Z' \cdot K_Y <0$.
\end{Lem}
\begin{proof}
Write $K_Y=g^*K_X+\sum_i a_i E_i$. Since $K_Y\cdot E_i\geq 0$, we have $a_i\leq 0$. It follows that $Z'\cdot K_Y=Z'\cdot g^*K_X+\sum_i a_iZ'\cdot E_i<0$.
\end{proof}

\begin{proof}[Proof of Proposition~\ref{cone2-1-1}] 
Again, we only prove the bound on $-(K_X+D)$ degree. Consider any $(K_X+D)$-negative extremal ray $R$. Choose a rational curve $C$ such that $R= \RR_{\geq 0}[C]$ and $C$ has the lowest $-(K_X+D)$-degree among all rational curves that generate $R$. We need to show  that $C \cdot (K_X+D) \geq -4$.  We prove this proposition by induction on the number of the irreducible components of $D=\sum _{i} a_i D_i $. 

When $D$ is empty, the statement is just Corollary~\ref{cone2-1}. 
Assume we have proved it when $D=\sum a_{i}D_{i}$ consists of  $k$ irreducible components.  
\newline \newline
Case I:  $C \cdot D_{i} \geq 0$ for some $i$.

We then have
\[C \cdot (K_{X}+D) = C \cdot (K_X+ \sum_j a_j D_j)\geq  C \cdot (K_{X}+ \sum_{j \neq i}a_{j}D_{j}).\] 

The log pair $(X, D'= \sum_{j \neq i}a_{j}D_{j})$ is also canonical. 
The inequality  $C \cdot (K_{X}+D) \geq -4$ follows by the induction hypothesis.
\newline \newline
Case II:  $C \cdot D_{i}<0$ for all $i$. 

Pick the divisor $D_1$, and write $B=\sum _{j >1}a_{j}D_{j}$. Note that $C \cdot D_1 <0$ by assumption. Therefore $C \subset D_1$. We prove by contradiction.
Assume now that $C \cdot (K_X+D) <-4$.
There are two subcases. 
\newline \newline
Case II-a: The pair $(X, D_{1}+B)$ is divisorial log terminal (dlt).

Note that $D_1$ is normal since $(X,D_1+B)$ is dlt. We have 
\[C \cdot (K_{X}+D) =C \cdot (K_X+a_1D_1+B) \geq  C \cdot (K_{X}+D_{1}+ B).\]
By adjunction, we have that $ C \cdot (K_{X}+D_{1}+B)= C \cdot (K_{D_{1}}+B^{D_1})$ where $B^{D_{1}}$ is the different of $B$ on $D_1$ (see \cite{ketc} for the definition of the different). 
The pair $(D_1, B^{D_1})$ is dlt, and the divisor $B^{D_1}$ is still a boundary divisor.

By our assumption, we have $C \cdot (K_X+D) < -4$. It follows that   $C \cdot (K_{D_{1}} +B^{D_{1}}) <-4 <-3$. By Lemma~\ref{sur1}, $C \sim \sum C_{i}$ 
(on $D_{1}$), all of $C_i$ are rational curves ($C_i$ may equal to $C_j$), and $C_{i} \cdot (K_{D_{1}} +B^{D_{1}}) \geq -3$. 
Pick one  $C_i$, say $C_1$. 
Since the ray $R=\RR_{\geq 0}[C]$ is extremal, $[C_1]$ also generates $R$. Therefore  $[C_1]$ is a positive multiple of $[C]$. It follows that $C_1 \cdot D_1 <0$, and 
\[ C_1 \cdot (K_X+D ) \geq C_1 \cdot (K_X+D_1+B)= C_1 \cdot (K_{D_1}+B^{D_1}) \geq -3.\] 
This contradicts to the choice of $C$. Therefore $$C\cdot (K_X+D)\geq C\cdot (K_{D_1}+B^{D_1})\geq -3,$$ again a contradiction.
\newline \newline
Case II-b:  The pair $(X, D_{1}+B)$ is not dlt. 

Again we start with the  assumption  $C \cdot (K_X+D) <-4$.
Consider a log terminal partial resolution $g:(Y,D'_1+B') \to (X,D_1+B)$ such that $(Y, D'_1+B'+F)$ is dlt, where $F=\sum_{j \in J}F_j$, $\{F_j: j \in J\}$ is the collection of all components of the exceptional divisor, and $K_Y+ D'_1+B'+F$ is $g$-nef (see Theorem 6.16 of \cite{ketc}). %%% check reference%% 
Write  
\[K_{Y}+D'_{1}+ B'+ \sum_{j \in J} F_j = g^{*}(K_{X}+D_1+B)+   \sum_{j \in J} b_{j}F_{j}.\] 
Note that  $b_j \leq 0$ since $K_{Y}+D'_{1}+ B'+ \sum _{j \in J}F_j$ is $g$-nef. 
 
We consider the following subcases:
\newline \newline
Case II-b-1:   $ C \nsubseteq \cup_{j \in J} g_{*}F_j$. 

Take the proper transform $C'$ of $C$ in $Y$. Note that $C' \cdot F_j \geq 0$ since $C' \nsubseteq F_j$, and 
\[ C' \cdot(K_Y+D'_1+B'+F)= C' \cdot g^*(K_X+D_1+B)+ C' \cdot (\sum _{j \in J}b_j F_j) \]
\[ \leq C' \cdot g^* (K_X+D_1+B)= C \cdot (K_X+D_1+B)\leq C\cdot (K_X+D).\]
It suffices to show that  $C' \cdot(K_Y+D'_1+B'+F) \geq -3$.  
As in the proof of Lemma~\ref{lc-1}, we write $[C']= \sum b_i [C_i] + u^+$ and 
$g_*([C_i])=e_i [C]$, $g_*(u^+)=e [C]$. Note that  
$C_i \cdot (K_Y+D'_1+B'+F) \geq -3$ by applying Case II-a to the log pair $(Y, D'_1+B' +\sum _{j \in J}F_j)$. Also, if $C_i$ is not contracted by $g$ then
$e_i\geq 1$ (otherwise $g_*(C_i)$ has smaller $-(K_X+D)$ degree).

Set $I_1=\{ i: C_i\;\; \text{is not contracted by $g$}\}$ and $I_2= \{ i: C_i\;\; \text{is  contracted by $g$}\}$. Recall that $(K_Y+D'_1+B'+F)$ is $g$-nef, therefore $C_i \cdot (K_Y+D'_1+B'+F) \geq 0$ when $C_i$ is contracted by $g$, i.e. $i \in I_2$. 
It is clear that  $\sum_{i \in I_1} b_i e_i+e=1$ since $g_* [C']=[C]$. We have  
$$ C' \cdot(K_Y+D'_1+B'+F)= [\sum_{i \in I_1}  b_i [C_i]+ \sum_{j \in I_2} b_j [C_j]  +u^+] \cdot(K_Y+D'_1+B'+F)$$
$$\geq (\sum_{i \in I_1}  b_i [C_i]) \cdot (K_Y+D'_1+B'+F)\geq -3 \sum_{i \in I_1} b_i.$$
The first inequality follows since $C_j$ is $g$-nef when $j \in I_2$.
Since $$\sum_{i \in I_1} e_i b_i +e=1,\;  e_i \geq 1, \; \text{and}\; e \geq 0,$$ it follows that 
$$\sum_{i \in I_1} b_i \leq 1 \; \text{ and hence}\; -3 \sum_{i \in I_1} b_i \geq -3.$$
This concludes this case.
\newline \newline 
Case II-b-2:  $C \subset g_{*}F_j$ for some $j \in J$. 

Consider the non-empty subset $J_0=\{j: C \subset g_{*}F_j\} \subset J$. 
Write \[K_Y+B'=g^{*}(K_X+B)+ \sum_{j \in J} k_j F_j, \quad g^{*}D_1=D'_1 + \sum_{j \in J} c_j F_j. \]
Consider $\mu =min_{j \in J_0}\{\frac{k_j+1}{c_j}\}$. Let $J_1 =\{ j \in J_0:  \frac{k_j+1}{c_j}=\mu\}$.
With this choice of $\mu$, the log pair $(X, \mu D_1+B)$ is strictly log canonical at a general point $p$ on $C$. Consider the $(K_Y+\mu D'_1+B + \sum_{j \in J} F_j)$-MMP over $(X, \mu D_1+B)$. Note that the pair   $(Y, \mu D'_1+B' + \sum_{j \in J} F_j)$ is dlt.
In each step of this MMP,  we either contract a divisor or perform a flip.     
Consider  the final outcome of this MMP \[g_1: (Y_1, \mu D'_1+B' + \sum _{j \in J_2} F_j) \to (X, \mu D_1+B),\] 
where $J_2$ is a subset of $J$ (that is, the set of all $j$ so that $F_j$ is not contracted in every step of the MMP). Note that since $X$ is $\QQ$-factorial, the exceptional locus of $g_1:Y_1 \to X$ is of codimension $1$. The pair $(Y_1, \mu D'_1+B' + \sum_{j \in J_2} F_j)$ is $\QQ$-factorial and dlt, and $K_{Y_{1}}+\mu D'_1+B' + \sum_{j \in J_2} F_j$ is $g_1$-nef. Observe that every irreducible component of the exceptional divisor on $Y_1$ is the  proper transform of some $F_j$. Also note that at least one of $F_j$ with  $j \in J_1$ (in fact, all of them) survives in every step of the MMP  since 
$(Y_1, \mu D'_1+B' +F)$ is dlt, and $(X, \mu D_1+B)$ is strictly log canonical at a general point $p$ on $C$. 

Write 
\begin{equation}\label{remif}
K_{Y_{1}} +\mu D'_1+B'+\sum _{j \in J_2} F_j =g_1^{*}(K_X+ \mu D_1+ B)+ \sum_{j \in J_2} f_j F_j.
\end{equation}
Since $K_{Y_{1}} +\mu D'_1+B'+\sum _{j \in J_2} F_j$ is $g_1$-nef, it follows that  $f_{j} \leq 0$.

Let $J_3=\{ j \in J_2: C \nsubseteqq  g_{1*} F_j \} \subset J_2$. Since $C \nsubseteqq \cup_{j \in J_3} g_{1*}  F_j$,  we can pick  a point $p \in C \setminus\{\cup_{j \in J_3}g_{1*} F_j\}$. 
Consider a small neighborhood $U \subset X$ of $p$ such that $U \cap \; ( \cup_{j \in J_3}g_{1*} F_j)=\varnothing$.

The pair $(X, \mu D_1+B)$ is log canonical on $U$. It follows that $f_j \geq 0$ when $C \subset g_{1*}F_j$, i.e. when $j \in J_2 \setminus J_3$.
 Therefore $f_j=0$ if $g_{*}F_j$ contains $C$. 

Since $f_j=0$ when $j \in J_2 \setminus J_3$ we can rewrite (\ref{remif}) as
%\[ K_{Y_1}+\mu D'_1+B'+\sum _{j \in J_2}F_j=g_1^*(K_X+\mu D_1+B)+\sum _{j \in J_2} f_jF_j\] as
\[ K_{Y_1}+\mu D'_1+B'+\sum _{j \in J_2}F_j=g_1^*(K_X+\mu D_1+B)+\sum _{j \in J_3} f_jF_j.\] 

Pick a divisor $F_j$ such that $C \subset g_{1*}F_j$. Take $Z \subset F_j \subset Y_1$ such that $g_{1*}Z=d C$. Note that $Z \cdot F_k \geq 0$ and $f_k \leq 0$ for all $k \in J_3$. It follows that  
\[ Z \cdot [K_{Y_1} +\mu D'_1+B'+\sum _{k \in J_2} F_k]= Z \cdot [g_1^*(K_X+\mu D_1+B)+\sum _{k \in J_3} f_k F_k] \] 
\[ \leq Z \cdot g_1^*(K_X+\mu D_1+B) \;\text{ (since $f_k \leq 0$ and $Z \cdot F_k \geq 0$ when $k \in J_3$)} \]
\[ = g_{1*} Z \cdot (K_X+\mu D_1+B) = d C \cdot (K_X+\mu D_1 +B)\] 
\[ \leq d C \cdot (K_X+ a_1 D_1 +B) \;\text{ (since $a_1 < \mu$ and $C \cdot D_1<0$})\]
\[= d C \cdot (K_X+ D)<-4.\] 

Our strategy is to break $Z$ in $F_j$. Note that $F_j$ is normal since $(Y_1, \mu D'_1+B'+ \sum_{ i \in J_2} F_i)$ is dlt. We need the next claim:
\begin{cl}\label{ratlsec}
The curve $Z$ can be taken to be a rational curve.
\end{cl}

Suppose that $Z$ is rational.
By adjunction, 
\[Z \cdot [K_{F_{j}}+ \mu (D'_1)^{F_j}+(B')^{F_j}+ (\sum_{i \in J_2, i \neq j}F_i)^{F_j}]\]
\[= Z \cdot [K_{Y_1}+ \mu (D'_1)+(B')+ \sum_{i \in J_2}F_i ] <-4.\]
Applying Lemma~\ref{sur1}, we have $Z \sim \sum Z_i$ with every $Z_i$ rational 
and \[Z_i \cdot (K_{F_j} +\mu (D'_1)^{F_j}+(B')^{F_j}+[\sum _{k \in J_2, \;k \neq j} F_k ]^{F_j}) \geq -3.\] 
Since $g_{1*}Z=dC$, it follows that there is at least one $Z_i$, say $Z_1$, which is not contracted to a point by $g_{1}\mid_{F_j}$. We can choose $Z_1$ such that $Z_1 \nsubseteqq F_l$ when $l \in J_3$.  
%since $C \nsubseteqq g_{1*}F_l$. 
Since $C$ represents an extremal ray, the class $g_{1*}[Z_1]$ is proportional to $[C]$ and $g_{1*}[Z_1]\cdot D_1<0$. Recall that $f_l \leq 0$ when $l \in J_3$. It follows that 
$$g_{1*}Z_1 \cdot (K_X+D)\geq g_{1*}Z_1\cdot (K_X+\mu D_1+B)$$
$$=Z_1 \cdot (K_{Y_{1}} +\mu D'_1+B'+\sum _{k \in J_2} F_k -\sum_{l \in J_3} f_l F_l ) \geq -3$$ since $Z_1 \cdot  F_l \geq 0$ and $-f_l \geq 0$. So the rational curve $g_{1*}Z_1$ has $-(K_X+D)$ degree $\leq 3$, contradicting to the assumptions that $C$ has the lowest $-(K_X+D)$-degree among the rational curves generating the ray $R$ and $C \cdot (K_X+D)<-4$.
\end{proof}

\begin{proof}[Proof of Claim \ref{ratlsec}]

Assume $Z$ is not rational. Consider the morphism $g_1: Y_1 \to X$ restricted to $F_j \subset Y \to C \subset X$. Note that the normalization of $C$ is $P^1$. The morphism $g_1\mid_{F_j}: F_j \to C$ factors through $F_j \to P^1 \to C$. Let $S \to F_j$ be the minimal resolution. We divide into two subcases:
\newline \newline
Subcase I: The general fiber of $S \to P^1$ is rational. 

The surface $S \to P^1$ is a ruled surface.
There is a section $\sigma: P^1 \to S$ with image $C_0$. Replace $Z$ by $C_0$ in $F_j$,
\newline \newline
Subcase II: The general fiber of $S \to P^1$ is irrational.

By Lemma~\ref{sur1-1}, we have   
$Z \cdot K_{F_j} <0$. Let $Z'$ be the proper transform of $Z$ in $S$.
It follows that $Z' \cdot K_S < 0$ by Lemma~\ref{sur2}. 
We can find a rational curve $C_x$ through a general point of $Z'$. Since the general fiber has no rational curve, the rational curve $C_x$ is not contained in any fiber, and is mapped onto $P^1$. Replace $Z$ by the image of $C_x$ in $F_j$. 
\end{proof}

\begin{Rem}
This induction process is quite complicated. If we apply the following theorem of Kawamata, the proof will be much easier.

\begin{Thm}[Kawamata \cite{ka91}]
Let $X$ be a normal projective variety with a boundary $\QQ$-divisor $D= \sum d_{i}D_i$, $0 \leq d_i \leq 1$, such that the log pair $(X,D)$ has only $\QQ$-factorial and log terminal singularities. Let $\phi: X \to Y$ be a morphism such that $-(K_X+D)$ is $\phi$-ample. Then every irreducible component $E$ of $Exc( \phi)$ is covered by rational curves $l$ such that 
\[ 0 < -(K_X+D) \cdot l \leq 2(dim\;E - dim \phi(E))\]
with $\phi(l)= {\text a \;point\; on\; Y}.$
\end{Thm}

Applying this to the contraction $\phi$ of an extremal ray $R$. The only case of Proposition~\ref{cone2-1-1} that is not covered by Kawamata's theorem is when $Exc(\phi)=X$ and $dim \phi(E)=0$. In this case, all curves in $X$ are contracted by $\phi$. So any curve represents the ray $R$. We may choose a curve $C$ not lying on $\cup_i D_i$. It follows that $C\cdot D_i\geq 0$ and $C \cdot (K_X+D) \geq C \cdot K_X$. We conclude by Corollary~\ref{cone2-1}.
\end{Rem}

{\bf Step 4:} The pair $(X,D)$ is klt.
\begin{Rem} 
Note that a dlt pair $(X,D)$ is the limit of klt pairs, that is, consider the limit of $(X, (1 -\epsilon) D)$ when $\epsilon \to 0$.    
The dlt case follows immediately from  the next corollary.
\end{Rem}
\begin{Cor}\label{cone2-2}
Let $(X,D)$ be a projective threefold pair, and $D$ a boundary  divisor. Assume the pair  $(X,D)$   has only Kawamata log terminal singularities. Then the set ${\bf \mR}$ of all $(K_{X}+D)$-negative extremal rays of $\overline{NE}(X)$ is countable and 
\[ \overline{NE}(X)=  \overline{NE}(X)_{K_{X}+D \geq 0} + \sum_{R_{i} \in  \mR}  R_i.\]
These rays are locally finite in the half-space $N_{1}(X)_{K_{X}+D <0}$, and for each ray $R_i$ there is a rational curve $l_i$ generating $R_i$ such that
\[l_i \cdot (-(K_{X}+D)) \leq 3+1.\]
\end{Cor} 

\begin{proof} We do induction on the number of negative discrepancies. 
Suppose $(X,D)$ is klt and for every exceptional divisor $E$, the discrepancy $a(E,X,D) \geq 0$, then it is canonical. The result follows from  Corollary~\ref{cone2-1}. 

Suppose $(X,D)$ has $k$ negative discrepancies. Let 
$$-d:=min\{a(E,X,D): E \text{ is an exceptional divisor on (a model of) } X \}<0.$$ 
Then we can find an extraction $g: (Y,g^{-1}_{*}(D)) \to (X,D)$ such that $K_{Y}+g^{-1}_{*}(D)+  d \sum E_{i}$ is $g$-nef, where $E_{i}$'s are all exceptional divisors (see \cite{ketc} Theorem6-16). Note that $K_{Y}+g^{-1}_{*}(D)+  d \sum E_{i}=g^{-1} (K_{X}+D)$ by our choice of $d$. 
Also observe that  $(Y, g^{-1}_{*}(D)+ d \sum E_{i})$ has fewer negative discrepancies. An argument similar to that in the proof of Lemma \ref{lc-1} completes the induction step.
\end{proof}

\begin{Rem}
This induction  does not work for  dlt or lc pairs since they may have infinitely many negative discrepancies. 
\end{Rem}


\begin{thebibliography}{HLOY02}
\bibitem[ACV03]{acv} D. Abramovich, A. Corti, and A. Vistoli, \emph{Twisted bundles and admissible covers}, Comm.\ Algebra 31 (2003) 3547--3618.

\bibitem[AV02]{av02} D. Abramovich, A. Vistoli. \emph{Compactifying the space of stable maps}, J. Amer. Math. Soc. 15, no. 1 27-75, 2002.

\bibitem[Am03]{am03} F. Ambro.  \emph{Quasi-log varieties}, Tr. Mat. Inst. Steklova 240 (2003), Biratsion. Geom. Linein. Sist. Konechno Porozhdennye Algebry, 220-239; translation in Proc. Steklov Inst. Math. 2003, no. 1(240), 214-233, 2003.

\bibitem[C04]{ca04} C. Cadman.   \emph{Counting Rational Curves with Perscribe Tangencies to a Smooth Plane Cubis}, preprint, 2004.  

\bibitem[De01]{de} O. Debarre.   \emph{Higher-dimensional algebraic geometry}, Springer-Verlag, 2001.  

\bibitem[FuP97]{FuP} W. Fulton and R. Pandharipande, \emph{Notes on stable maps and quantum cohomology}, in {\sl Algebraic geometry (Santa Cruz, 1995)}, 45--96, Amer.\ Math.\ Soc., 1997.

\bibitem[Gi84]{gi84} H. Gillet, \emph{Intersection theory on algebraic stacks and $Q$-varieties}, J. Pure Appl. Algebra 34 (1984), no. 2-3, 193--240.

\bibitem[I71]{i71} L. Illusie, \emph{Complexe cotangent et d\'eformations I}, Lecture Notes in Mathematics, Vol. 239, Springer-Verlag, 1971.

\bibitem[Ka84]{ka84}
Y. Kawamata, \emph{The cone of curves of algebraic varieties}, Ann. of Math. (2) 119 (1984), no. 3, 603--633.

\bibitem[Ka91]{ka91} Y. Kawamata, \emph{On the length of an extremal rational curve}, Invent. Math. 105 (1991), no. 3, 609--611. 

\bibitem[Ka95]{ka95} Y. Kawamata. \emph{Unobstructed deformations II}, J. Algebraic Geometry 4 (1995) 277-279.

\bibitem[Kw79]{Kw79} Y. Kawasaki. \emph{The Riemann-Roch theorem for complex V-manifolds}, Osaka J. Math.16 (1979) 151-159.

\bibitem[Ke99]{ke99} S. Keel. \emph{Basepoint freeness for big line bundles in positive characteristic}, preprint math.AG/9901149, Ann. of Math. (2) 149 (1999), no. 1, 253-286, 1999. 

\bibitem[Ko84]{ko84} J. Koll\'ar. \emph{The cone theorem}, Ann. of Math. (2) 120 (1984),  no. 1, 1--5. 

\bibitem[Ko91]{ko91} J. Koll\'ar. \emph{Flips, flops, minimal models}, etc, Surv. in Diff. Geom., 1:113-199, 1991.

\bibitem[K et al92]{ketc} J. Koll\'ar et al, \emph{Flips and Abundance for Algebraic Threefolds}, vol. 211, Ast\'erisque, 1992. 

\bibitem[Ko92]{ko92} J. Koll\'ar. \emph{Cone Theorems and Bug-eyed covers}, J. Algebraic Geometry, 1:293-323, 1992.

\bibitem[Ko96]{ko96} J. Koll\'ar. \emph{Rational curves on algebraic varieties}, Springer-Verlag, 1996.  

\bibitem[KM98]{km98} J. Koll\'ar, S. Mori. \emph{Birational Geometry of Algebraic Varieties}, Cambridge Tracts in Mathematics 134, Cambridge University Press, 1998.

\bibitem[LMB00]{LMB} G. Laumon and L. Moret-Bailly, \emph{Champs algebriques}, Springer-Verlag, 2000.

\bibitem[Ma02]{ma02} K. Matsuki, \emph{Introduction to the Mori program}, Springer-Verlag, 2002.

\bibitem[MO02]{mo02} K. Matsuki, M. Olsson. \emph{Kawamata-Viehweg vanishing as Kodaira vanishing for stacks}, preprint math.AG/0212259, to appear in Math. Res. Letters.

\bibitem[Mo82]{mo82} S. Mori. \emph{Threefolds whose canonical bundles are not numerically effective}, Ann. of Math. (2) 116 (1982), no. 1, 133--176.

\bibitem[O03]{o03} M. Olsson, \emph{Hom--stacks and restriction of scalars}, preprint, 2003.

\bibitem[R83]{r83} M. Reid, \emph{Projective morphisms according to Kawamata}.

\bibitem[Sh85]{sh85} V. V. Shokurov, \emph{A nonvanishing theorem}, Izv. Akad. Nauk SSSR Ser. Mat. 49 (1985), no. 3, 635--651.

\bibitem[T99]{t99} B. Toen, \emph{Th\'eor\`emes de Riemann-Roch pour les champs de Deligne-Mumford}, $K$-Theory 18 (1999), no. 1, 33--76. 

\bibitem[To02]{to02} B. Totaro. \emph{The resolution property for schemes and stacks}, preprint math.AG/0207210, J. Reine Angew. Math. 577 (2004), 1--22.

\bibitem[Vi89]{vi89} A. Vistoli, \emph{Intersection theory on algebraic stacks and on their moduli spaces} Invent. Math. 97 (1989), no. 3, 613--670. 

\end{thebibliography}
\end{document}